\documentclass[10pt]{amsart}
\usepackage{amsfonts,latexsym, color}
\usepackage{amsmath}
\usepackage{amssymb}
 \font\tenmsb=msbm10 at 12pt \font\sevenmsb=msbm7 at 8pt \font\fivemsb=msbm5 at
6pt
\newfam\msbfam
\textfont\msbfam=\tenmsb \scriptfont\msbfam=\sevenmsb \scriptscriptfont\msbfam=\fivemsb

\def\R{{\mathbb R}}

\def\N{{\mathbb N}}
\def\Sp{{\mathbb S}}

\begin{document}
\newcommand{\reset}{\setcounter{equation}{0}}

\newcommand{\beq}{\begin{equation}}
\newcommand{\noi}{\noindent}
\newcommand{\eeq}{\end{equation}}
\newcommand{\dis}{\displaystyle}
\newcommand{\mint}{-\!\!\!\!\!\!\int}

\def \theequation{\arabic{section}.\arabic{equation}}

\newtheorem{thm}{Theorem}[section]
\newtheorem{lem}[thm]{Lemma}
\newtheorem{cor}[thm]{Corollary}
\newtheorem{prop}[thm]{Proposition}
\theoremstyle{definition}
\newtheorem{defn}{Definition}
\newtheorem{rem}[thm]{Remark}

\def \bx{\hspace{2.5mm}\rule{2.5mm}{2.5mm}} \def \vs{\vspace*{0.2cm}} \def
\hs{\hspace*{0.6cm}}
\def \ds{\displaystyle}
\def \p{\partial}
\def \O{\Omega}
\def \b{\beta}
\def \m{\mu}
\def \M{{\mathcal M}}
\def \G{{\mathbb G}}
\def \T{{\mathbb T}}
\def \ou{{\overline u}}
\def \ov{{\overline v}}
\def \D{\Delta}
\def \dD{{\mathcal D}}
\def\H{{\mathbb H}}
\def \cR {{\mathcal R}}
\def \d{\delta}
\def \s{\sigma}
\def \e{\varepsilon}
\def \a{\alpha}
\def \L{{\mathcal L}}
\def \o{\omega}
\def \g{\gamma}
\def \wv{{\widetilde v}}
\def \Sp{{\mathbb S}}
\def \hw{{\widehat w}}
\def\cqfd{%
\mbox{ }%
\nolinebreak%
\hfill%
\rule{2mm} {2mm}%
\medbreak%
\par%
}
\def \pr {\noindent {\it Proof:} }
\def \rmk {\noindent {\it Remark} }
\def \esp {\hspace{4mm}}
\def \dsp {\hspace{2mm}}
\def \ssp {\hspace{1mm}}

\title{Higher order Hardy-Rellich identities}
\author{Xia Huang}
\address{School of Mathematical Sciences,  Key Laboratory of MEA (Ministry of Education) \& Shanghai Key Laboratory of PMMP,  East China Normal University, Shanghai 200241, China}
\email{xhuang@cpde.ecnu.edu.cn}
\author{Dong Ye}
\address{School of Mathematical Sciences,  Key Laboratory of MEA (Ministry of Education) \& Shanghai Key Laboratory of PMMP,  East China Normal University, Shanghai 200241, China}
%\address{IECL, UMR 7502, University of Lorraine, 57050 Metz, France}
\email{dye@math.ecnu.edu.cn}

\date{}
\begin{abstract}
In this paper, we show Hardy-Rellich identities for polyharmonic operators $\Delta^m$ and radial Laplacian $\Delta_r^m$ in $\R^n$ with Hardy-H\'enon weight $|x|^\alpha$ for all $m, n\in \N, \alpha\in \R$. Moreover, the iterative method is applied to give Hardy-Rellich equalities with general weights on Riemannian manifolds. These identities provide naturally an alternative approach to obtain and improve Hardy-Rellich type inequalities. As example of application, we extend several Rellich inequalities of Tertikas-Zographopoulos \cite{TZ} to the weighted case; using equality with weights involving logarithmic, we show another new weighted Rellich estimate between integrals of $\Delta u$ and $|\nabla u|$; we establish also a Rellich identity involving the Laplace-Beltrami operator $\Delta_\H$ and the radial Laplacian $\Delta_{\rho, \H}$ of the hyperbolic space $\H^n$, which yields in particular
%$\|\Delta_\H u\| \geq \|\Delta_{\rho,\H} u\|$ for any $u \in C_c^2(\H^n\backslash\{0\})$ with $n \geq 3$, and
brand-new Rellich inequalities for $\|\Delta_\H u\|$ in $\H^3$ and $\H^4$.
\end{abstract}

\maketitle
\noindent
{\small {\bf 2010 MSC:} {\sl 26D10, 26D15, 35A23}}

\noindent
{\small {\bf Keywords:} {\sl higher order Hardy-Rellich identities, polyharmonic operators, iteration}}

\vskip 0.8cm

\section{Introduction}
The first order Hardy inequalities go back to G.H. Hardy, who showed in \cite{GH} a very famous estimate. Let $p>1$, then
$$
\int_{0}^\infty |u'(x)|^p dx \geq \left(\frac{p-1}{p}\right)^p \int_0^\infty \frac{|u(x)|^p}{x^p} dx, \quad \forall~u\in C^1(\mathbb{R}_+),~u(0)=0.
$$
Since one century, Hardy type inequalities have been enriched extensively and broadly, they play an important roles in many branches of analysis and geometry. A huge literature exists, we refer to the classical or more recent books \cite{BEL, GM1, Ma, OK, RS1} for interested readers.

\medskip
Now we have various methods to provide Hardy type inequalities. One way is to establish {\it equalities} with nonnegative remainder terms which yield immediately the corresponding inequalities. Look at the most famous weighted Hardy inequality:
\begin{align}
\label{cHR1a}
\int_{\R^n} \frac{|\nabla u|^2}{|x|^\alpha} dx \geq \frac{(n-2 - \alpha)^2}{4}\int_{\R^n}\frac{u^2}{|x|^{\alpha +2}} dx, \quad \forall \; \alpha \in \R, \; u \in C^1_c(\R^n\backslash\{0\}).
\end{align}
Indeed, \eqref{cHR1a} is a direct consequence of the following identity: For $n \ge 1$, $u\in C^1_c(\R^n\backslash\{0\})$ and $\a \in \R$, there holds
\begin{align}
\label{HR1a}
\|\nabla u\|_\a^2 - \frac{(n-2 - \alpha)^2}{4}\|u\|_{\a+2}^2 = \left\|\nabla\left(r^\frac{n-2-\alpha}{2}u\right)\right\|_{n-2}^2.
\end{align}
Here $r = |x|$ stands for the Euclidean norm in $\R^n$, and
\begin{align}
\label{weight1}
\|g\|_\b^2 := \int_{\R^n} \frac{|g(x)|^2}{r^\beta} dx, \quad \forall\; g \in C_c(\R^n\backslash\{0\}), \; \b \in \R.
\end{align}
The equality \eqref{HR1a} can be found for example in \cite{DV} with $\alpha = 0$; or \cite[Lemma 2.3(i)]{TZ} with $\alpha = 2$. On the other hand, we can also use radial derivative to derive \eqref{cHR1a}, since
\begin{align}
\label{rHR1a}
\|\p_r u\|_\a^2 - \left(\frac{n-2 -\alpha}{2}\right)^2\|u\|_{\a+2}^2 = \left\|T_\a(u)\right\|_\a^2, \quad \forall\; u \in C^1_c(\R^n\backslash\{0\}), \; \a \in \R,
\end{align}
where
\begin{align}
\label{TaLj}T_\a(f) := \p_r f + \frac{n-2-\alpha}{2r}f \quad \mbox{and} \quad \p_r = \frac{x}{r}\cdot\nabla = \sum_{i= 1}^n \frac{x_i}{r}\p_{x_i}.
\end{align}
Both \eqref{HR1a} and \eqref{rHR1a} can be proved easily using integration by parts (see Remark \ref{newrmk1} below). The formula \eqref{rHR1a} was found firstly by Brezis-V\'azquez with $\alpha=0$, they called it the {\it magical computation}, see \cite[page 454]{BV}.

\smallskip
In 2011, Ghoussoub-Moradifam introduced the powerful idea of Bessel pairs in \cite{GM} (see also \cite{OK}) to establish Hardy type inequalities
$$\int_\Omega V|\nabla u|^2 dx \geq \int_\O Wu^2dx,$$
with more general {\it radial} weights $V$ and $W$. Their technique reaches a large success, and has been applied in many works to prove and improve first order Hardy inequalities, to name a few, see for instance \cite{BGR, CFLL, DLT3, FLL, FLLM, GR} and references therein.

\smallskip
Later on, we attempt a more general approach in \cite{HY}, where the departure point is an elementary equality as follows. Let $(\mathcal{M}, g)$ be a Riemannian manifold. Consider $V \in C^1(\M)$ and ~$\vec{F}\in C^1(\M, T_g \mathcal{M})$, then for any $u\in C^1_c(\M)$, there holds
\begin{align*}
\int_\M V\|\nabla_g u\|^2 dg = \int_\M \big[{\rm div}_g (V\vec{F}) -V\|\vec{F}\|^2\big] u^2 dg +\int_\M V\|\nabla_g u + u\vec{F}\|^2 dg.
\end{align*}
In particular, let $\vec{F}=-\frac{\nabla_g f}{f}$ with $f$ positive in $C^2(\M)$, we see that
\begin{align}
\label{HYequa1}
\int_\M V\|\nabla_g u\|^2 dg = -\int_\M \frac{{\rm div}_g (V\nabla_g f)} {f} u^2 dg +\int_\M V\Big\|\nabla_g u -\frac{u}{f}\nabla_g f\Big\|^2 dg.
\end{align}
We deduce then, for $V \geq 0$,
\begin{align*}
\int_\M V\|\nabla_g u\|^2 dg \ge -\int_\M \frac{{\rm div}_g (V\nabla_g f)} {f} u^2 dg, \quad \forall\; u\in C^1_c(\M).
\end{align*}
Hence, by choosing suitable auxiliary function $f$, we can refind, extend or improve most of first order Hardy inequalities by the identity \eqref{HYequa1}, see \cite{HY} for many examples of application.
%see \cite{HY} and \cite{CZ}, \cite{FTT}. %\footnote{We should add some references of other authors}.

\smallskip
In particular, when $(\M, g) \subset \R^n$, $V$ and $f$ are positive and radial, we find the radial Bessel pair of Ghoussoub-Moradifam with $V$ and
$$W =  - \frac{{\rm div}(V\nabla f)} {f} =  -\frac{V}{f}\Big[f'' + \Big(\frac{n-1}{r} + \frac{V'}{V}\Big)f'\Big].$$

It is also worthy to mention that Kaj\'ant\'o-Krist\'aly-Peter-Zhao use recently the Riccati pairs to extend the Bessel pairs, and provide an alternative approach to get Hardy inequalities on complete, non-compact Riemannian manifolds, see \cite{KKPZ}.

\medskip
The well-known second order Hardy-Rellich inequality in $\R^n$ is the following (see \cite[Theorem 4.1]{CM}): Let $n\ge 2$, $|\a +2|\leq 2\sqrt{n^2 -2n+2}$, $n\ge 2$,
\begin{align}
\label{cHR2a}
\int_{\R^n} \frac{|\D u|^2}{|x|^\alpha} dx \geq \left[\frac{(n + \alpha)(n-4 - \alpha)}{4}\right]^2\int_{\R^n}\frac{u^2}{|x|^{4 + \alpha}} dx, \quad \forall\; u \in C_c^2(\R^n\backslash\{0\}).
\end{align}
The inequality \eqref{cHR2a} with $\alpha=0$, $n\ge 5$ was first given by Rellich at the ICM 1954 \cite{Re}, that's the reason that Hardy type estimates with higher order differential operators are named often Hardy-Rellich inequalities. In 2017, Machihara-Ozawa-Wadade reproved \eqref{cHR2a} with $\alpha=0$ elegantly using two identities, see \cite[Theorems 1.1-1.2]{MOW}. The idea in \cite{MOW} works indeed for general weights $|x|^\a$ with $\a \in \R$ (see section 3 and also \cite[Theorem 4.1]{RS} for \eqref{rHR2a}). The first equality shows the relationship between $\|\D u\|_\a$ and $\|\D_r u\|_\a$,
\begin{align}
\label{HR2a}
\begin{split}
\|\D u\|_\a^2 & = \|\D_r u\|_\a^2 + \left\|(\Delta - \D_r) u\right\|_\a^2\\ &\quad + \frac{(n + \alpha)(n-4 - \alpha)}{2}\sum_{j=1}^n\left\|L_ju\right\|_{\a+2}^2 + 2\sum_{j=1}^n\left\|T_\a\left(L_ju\right)\right\|_\a^2.
\end{split}
\end{align}
Here $$\Delta_r = \p_r^2 + \frac{n-1}{r}\p_r,$$
is the radial Laplacian in $\R^n$, and
\begin{align}
\label{defLj}
L_j := \p_{x_j}  - \frac{x_j}{r}\p_r, \;\;1\leq j \leq n.
\end{align}
The second equality is
\begin{align}
\label{rHR2a}
\begin{split}
\|\D_r u\|_\a^2 = & \; \frac{(n+\a)^2(n-4-\a)^2}{16}\left\|u\right\|_{\a+4}^2 \\
& + \left[\frac{(n+\a)^2}{4} + \frac{(n-4-\a)^2}{4}\right]\left\|T_{\a+2}(u)\right\|_{\a+2}^2 + \left\|T_\a\circ T_{\a+2}(u)\right\|_{\a}^2.
\end{split}
\end{align}
The paper of Machihara-Ozawa-Wadade is the first work which provides second order Rellich type equalities, their idea is generalized for various situations such as for more general {\it radial} weights in $\R^n$, in half space $\R^n_+$, homogeneous groups or hyperbolic spaces, see \cite{BGR, BMO, DLT1, DLT3, DLP} and references therein.

\medskip
Another interesting technique using identity to get Rellich inequalities is the factorization method (see Gesztesy-Littlejohn \cite{GL}), where they decompose cleverly the positive operator $T_{\a, \b}^*T_{\a, \b}$ for $T_{\a, \b} = -\Delta + \a r^{-2}x\cdot\nabla + \b r^{-2}$ with $\a, \b \in \R$. But this approach seems to be difficultly adaptable for higher order situations.

\medskip
An interesting question arises: Can we find Hardy-Rellich equalities for general higher order operators?

\smallskip
The first notable example is due to Ruzhansky-Surugan \cite[(6.2) and Theorem 6.1]{RS}. They considered homogeneous Carnot group $\G$ and established a remarkable equality for $\|\cR^k \|_\a^2$ for all $k \in \N$, $\a \in \R$. Here $\cR =\frac{d}{d\|x\|}$ is the radial operator (with respect to a quasi-norm $\|\cdot\|$ over $\G$), a natural analogue to $\p_r$ in Euclidean spaces, and $\cR^k$ means the $k$-times composition of $\cR$.

\medskip
However, to our best knowledge, there is no general Hardy-Rellich type identities for polyharmonic operators $\Delta^m$ in the literature. Our first objective here is to show that such equalities exist.

\medskip
Let us begin with the radial polyharmonic Laplacian $\Delta_r^m$ in $\R^n$. To write properly higher order Hardy-Rellich identities, we denote
\begin{align}
\label{Rak}
R_{\a,k}:=T_\a \circ T_{\a+2}\circ\cdot\cdot\cdot\circ T_{\a+2k}, \quad \forall\; \a \in \R,\; k \in \N,
\end{align}
and
\begin{align}
\label{AD}
A_\beta=\frac{(n+\beta)(n-4-\beta)}{4},\;\; D_\beta=\frac{(n+\beta)^2+(n-4-\beta)^2}{4},\quad \forall~\beta\in\mathbb{R}.
\end{align}

\begin{thm}
\label{thmrHRma}
Let $m, n \geq 1$, $\a \in \R$ and $u \in C_c^{2m}(\R^n\backslash\{0\})$, there holds
\begin{align}
\label{rHRma}
\begin{split}
\|\D_r^m u\|_\a^2 = \prod_{\ell = 0}^{m-1} A_{\a+4\ell}^2 \times \|u\|_{\a+4m}^2 + \sum_{j =0}^{2m-1} C_{j,m,\a} \left\|R_{\a+2j,2m-1-j}(u)\right\|_{\a+2j}^2,
\end{split}
\end{align}
where the positive coefficients $C_{j,m,\a}$ $(0 \leq j \leq 2m-1)$ are given by the iterative formula
\begin{align}
\label{iteC}
C_{j, k+1, \a} = A_{\a+4k}^2C_{j-2, k, \a} + D_{\a + 4k}C_{j-1, k, \a} + C_{j, k, \a}, \quad \forall\; 0\leq j \leq 2k+2, k \geq 0,
\end{align}
with the convention $C_{0, 0, \a} = 1$, and $C_{j, k, \a} = 0$ for $j < 0$ or $j > 2k$.
\end{thm}

As all coefficients $C_{j, m, \a}$ in \eqref{rHRma} are positive, we get obviously the following Hardy-Rellich inequality with optimal coefficient, and see that the equality holds only if $u \equiv 0$. Let $\a \in \R$, $m\in \N$,
\begin{align*}
\|\D_r^m u\|_\a^2 \ge \prod_{\ell = 0}^{m-1} A_{\a+4\ell}^2 \|u\|_{\a+4m}^2, \quad \forall\;  u \in C_c^{2m}(\R^n\backslash\{0\}).
\end{align*}
Moreover, seeing \eqref{rHRma} and definition of $R_{\a, k}$, the remainder terms can be seen coming from lower order Hardy-Rellich inequalities. For example, look at the second order case $m=1$, we can rewrite \eqref{rHR2a} as
\begin{align*}
& \quad \|\D_r u\|_\a^2 - A_\a^2\|u\|_{\a+4}^2 \\
& = D_\a\|T_{\a+2}(u)\|_{\a +2}^2 + \|T_\a\circ T_{\a+2}(u)\|_\a^2\\
& = D_\a\left[\|\p_r u\|_{\a+2}^2 - \frac{(n-4 - \alpha)^2}{4}\|u\|_{\a+4}^2\right] + \left[\|\p_r(T_{\a+2}u)\|_\a^2 - \frac{(n-2- \alpha)^2}{4}\|T_{\a+2}u\|_{\a+2}^2\right],
\end{align*}
where both above bracket are remainder terms of Hardy inequality, given by \eqref{rHR1a}.

\medskip
On the other hand, we see that
\begin{align*}
%\label{Cma}
C_{0, m+1, \a} = 1\;\;\mbox{and}\;\; C_{2m, m, \a} = \prod_{\ell =0}^{m-1} A_{\a+4\ell}^2, \quad \forall\; m \geq 0.
\end{align*}
So \eqref{rHRma} can be rewritten as
\begin{align*}\|\D_r^m u\|_\a^2 = C_{2m, m, \a} \|u\|_{\a+4m}^2+\sum_{j =0}^{2m-1} C_{j,m,\a} \left\|R_{\a+2j,2m-1-j}(u)\right\|_{\a+2j}^2, \quad \forall\; u \in C_c^{2m}(\R^{2m}\backslash\{0\}).\end{align*}

\begin{rem}
Directly or using density argument, all involved equalities or inequalities in this paper could hold true for more general functions. But for simplicity, we will always work with smooth functions compactly supported out of eventual singularities of the weights.
\end{rem}

Similar identities for polyharmonic operators $\Delta^m $ can be established, although a little bit more involved.
\begin{thm}
\label{newthm1}
For any $n, m \geq 1$, $\a \in \R$ and $u \in C_c^{2m}(\R^n\backslash\{0\})$, there holds
\begin{align}
\label{HRm}
\begin{split}
& \quad \|\D^m u\|_\a^2 -  H_{m-1, \a} \|u\|_{\a+4m}^2\\
 & = \sum_{k =1}^m H_{k-2, \a}\Big[ D_{\a+4k-4}\left\|T_{\a+4k-2}(\D^{m-k}f)\right\|_{\a+4k-2}^2 + \left\|R_{\a+4k, 1}(\D^{m-k} u)\right\|_{\a + 4k}^2\Big]\\
&\quad + \sum_{k = 1}^m H_{k-2, \a} \Big[\left\|(\Delta - \D_r) (\D^{m-k} u)\right\|_{\a+4k}^2  + 2\sum_{j=1}^n \left\|T_{\a+4k}(L_j (\D^{m-k} u)\right\|_{\a+4k-2}^2\Big]\\
&\quad + 2\sum_{j=1}^n \Big(\sum_{k = 1}^{m}A_{\a+4k - 4}H_{k-2, \a}\Big)\left\|L_j(\D^{m-k} u)\right\|_{\a+4k-2}^2,
\end{split}
\end{align}
where $H_{k, \a} = C_{2(k+1), k+1, \a}$ for $k \geq -1$, i.e.~$H_{-1, \a} = 1$ and $H_{k, \a} = \prod_{0 \le \ell \le k} A_{\a+4\ell}^2$ for $k \in \N$.
\end{thm}

Now we want to generalize \eqref{HYequa1} to higher order cases, in other words, to establish Hardy-Rellich identities for polyharmonic operators on manifolds and general weights. Let $(\M, g)$ be a Riemannian manifold, denote
\begin{align*}
\dD_{2m} = \Delta_g^m, \;\; \dD_{2m+1} = \nabla_g \Delta_g^m, \quad  \forall\; m\in \N.
\end{align*}
Therefore the corresponding adjoint operators are respectively $\dD^*_{2m} = \dD_{2m}$ and $\dD_{2m+1}^* = -\Delta_g^m({\rm div}_g)$.

\begin{thm}
\label{thmHRV}
Let $k \geq 1$, $V\in C^k(\M)$. Assume that $f\in C^k(\M)$ satisfies $|\dD_{2\ell} f|>0$ in $\M$ for all $\ell \ge 0$, $2\ell < k$. Then for any $u\in C_c^k(\M)$, we have
\begin{align}
\label{HRV}
\begin{split}
&\quad \int_{\M} V\|\dD_k u\|^2 dg-\int_{\M}\frac{\dD_k^* (V\dD_k f)}{f}u^2 dg\\& =-2\sum_{\ell \ge 0, 2\ell\leq k-2}\int_{\M}\frac{\dD_{k-2-2\ell}^*(V\dD_k f)}{\dD_{2\ell} f}\Big\|\dD_{2\ell+1} u -\frac{\dD_{2\ell}u}{\dD_{2\ell} f}\dD_{2\ell+1}f \Big\|^2 dg \\ & \quad +\sum_{\ell \geq 1, 2\ell \le k-1}\int_{\M}\frac{\dD_{k-2\ell}^*(V\dD_k f)}{\dD_{2\ell} f}
\Big|\dD_{2\ell} u -\frac{\dD_{2\ell}u}{\dD_{2\ell-2} f}\dD_{2\ell}f\Big|^2 dg\\
& \quad + \int_\M V\Big\|\dD_k u - \frac{\dD_{2[(k-1)/2]}u}{\dD_{2[(k-1)/2]}f}\dD_k f\Big\|^2 dg,
\end{split}
\end{align}
where $[\beta]$ denotes the integer part of $\beta \in \R$.
\end{thm}

A direct consequence is the following abstract Hardy-Rellich inequality.
\begin{cor}
\label{coriHRV}
Let $k\geq 1$ and $f\in C^k(\M)$ satisfy $(-\Delta_g)^\ell f > 0$ for $0 \le \ell \le k-1$, then
\begin{align}
\label{iHRV2k}
\int_{\M} \|\dD_k u\|^2 dx \geq \int_\M \frac{(-\Delta_g)^k f}{f} u^2 dx,\quad\forall~u\in C_c^k(\M).
\end{align}
Moreover, the equality holds if and only if $u \equiv 0$.
\end{cor}

Seeing \eqref{iHRV2k}, if $f$ presents suitable behavior (at infinity or/and near $\p\M$), we can expect to use cut-off function approximating $f$ and show eventually the optimality of some involved coefficients. For a very simple example, let $\M = \R^n\backslash\{0\}$ and $f=|x|^{\frac{2m-n}{2}}$ with $n > 4m$, we see that $(-\Delta)^\ell f > 0$ in $\R^n\backslash\{0\}$ for all $0 \le \ell \le 2m-1$. Hence by \eqref{iHRV2k}, we see the classical Hardy-Rellich inequality
\begin{align*}
\int_{\mathbb{R}^n} |\Delta^m u|^2 dx \geq H_{m, 0}\int_{\mathbb{R}^n} \frac{u^2}{|x|^{2m}} dx,\quad \forall~u\in C_c^{2m}(\mathbb{R}^n\setminus{\{0}\}).
\end{align*}
Of course, the above inequality is also a consequence of \eqref{HRm}, however with different expressions of (nonnegative) remainder terms. This means somehow the complexity and rich possibility of Rellich type inequalities, we need to pay attention to various remainder terms, if we aim to get different or improved estimates.

\medskip
As application, although we are convinced that much more results can be derived, we just expose several interesting consequences with second order differential operators. We will see that the Hardy-Rellich identities provide naturally new estimates, even with the Laplacian in $\R^n$. For example, the iterative approach to obtain Hardy-Rellich identities \eqref{rHRma} and \eqref{HRm} works also for weights involving logarithmic (see Proposition \ref{prop5.4} below), which yields
\begin{thm}
\label{thmHRln}
Let $n \ge 2$ and $u \in C_c^2(B_1\backslash\{0\})$ where $B_1$ is the unit ball in $\R^n$.
%If $|\alpha + 2| \leq \sqrt{n^2 -2n+2}$, we have
%\begin{align}\label{HR2ln}\|\D u\|_\a^2 &\geq A_\a^2 \left\| u\right\|^2_{\a+4} + D_\a\|(\ln|x|)^{-1}u\|_{\a+4}^2 + \frac{9}{16} \|(\ln|x|)^{-2} u\|_{\a+4}^2.\end{align}
If $|3\a + 4 +n| \leq 2\sqrt{n^2-n+1}$, there holds
\begin{align}
\label{HR21ln}
\|\D u\|_\a^2 \geq \frac{(n+\a)^2}{4} \left\|\nabla u\right\|^2_{\a+2} + \frac{(n-4-\a)^2}{16}\|(\ln|x|)^{-1}u\|_{\a+4}^2 + \frac{9}{16} \|(\ln|x|)^{-2} u\|_{\a+4}^2.
\end{align}
\end{thm}
Notice that by the study of Beckner \cite[Theorem 4]{B} and Tertikas-Zographopoulos \cite[Theorem1.7]{TZ} (see also the comments in \cite[page 83]{Hama}), the above constraint on $\alpha$ for \eqref{HR21ln} is the sharp range to claim
%$$\|\D u\|_\a^2 \geq \frac{(n+\a)^2(n-4-\alpha)^2}{16} \left\| u\right\|^2_{\a+4} \quad \mbox{or} \quad
\begin{align}
\label{HR21Rn}
\|\D u\|_\a^2 \geq \frac{(n+\a)^2}{4} \left\|\nabla u\right\|^2_{\a+2}, \quad \forall\; u \in C_c^2(B_1\backslash\{0\}).
\end{align}

As example on manifold, we consider the hyperbolic space $\H^n$, where the Rellich type inequalities have been studied intensively, see \cite{BGG, BGGP, BGR, BGR2, Ngu, KO, KO2} and references therein. In the spirit of Machihara-Ozawa-Wadade's formula for $\R^n$, we show a Rellich type identity on $\H^n$, which enable us quickly new Rellich inequalities.

\smallskip
More precisely, let $\H^n$ be Poincar\'e's hyperbolic ball and $\rho=\ln\frac{1+r}{1-r}$ be the geodesic distance from origin to $x\in \H^n$. Define $\L_j$ to be the sphere derivatives in $\H^n$,
\begin{align}
\label{Ljhyper}
\L_j=\frac{1-r^2}{2}\partial_{x_j}-\frac{x_j}{r}\partial_{\rho}=\frac{1-r^2}{2}L_j
\end{align}
where $L_j$ is that given in \eqref{defLj}. We denote by $\Delta_{\rho, \H}$ the radial Laplacian in $\H^n$, namely
\begin{align}
\label{Drhyper}
\Delta_{\rho, \H} = \frac{\partial^2}{\partial \rho^2} + (n-1)\coth \rho \frac{\partial}{\partial \rho}.
\end{align}
%We claim the following Rellich identity and inequalities.
\begin{thm}
\label{thm6}
For $n\geq 2$ and any $u \in C_c^2(\H^n\backslash\{0\})$, there holds
\begin{align}
\label{HR2=hyper}
\begin{split}
\|\Delta_\H u\|^2 & = \|\Delta_{\rho,\H} u\|^2 + \|(\Delta_\H -\Delta_{\rho,\H}) u\|^2\\
& \quad + 2\sum_{j=1}^n \|\partial_\rho (\L_j u)\|^2 -2\sum_{j=1}^n \|\L_j u\|^2 -2\sum_{j=1}^n\Big\|\frac{\L_j u}{{\rm sh}\rho}\Big\|^2.
\end{split}
\end{align}
Here $\|\cdot\|$ means the norm in $L^2(\H^n)$.
%\end{thm}
%\begin{thm}
%\label{compahyper}
%Let $n\geq 2$, $u \in C_c^2(\H^n\backslash\{0\})$,
Moreover, we have
\begin{align}
\label{HRhyper1}
\begin{split}
\|\Delta_\H u\|^2& \geq \|\Delta_{\rho,\H} u\|^2 +\frac{1}{2}\sum_{j=1}^n\Big\| \frac{\L_j u}{\rho}\Big\|^2 +\frac{(n+1)(n-3)}{2}\sum_{j=1}^n\|{\rm coth}\rho\L_j u\|^2
%+\frac{(n+1)(n-3)}{2}\sum_{j=1}^n\Big\|\frac{\L_j u}{{\rm sh}\rho}\Big\|^2,\\
\end{split}
\end{align}
Consequently $\|\Delta_\H u\| \geq \|\Delta_{\rho,\H} u\|$ for any $u \in C_c^2(\H^n\backslash\{0\})$ with $n \geq 3$.
%and the equality holds if and only if $u$ is radial.
\end{thm}

Notice that the comparison result $\|\Delta_\H u\| \geq \|\Delta_{\rho,\H} u\|$ in $C_c^2(\H^n\backslash\{0\})$ was proved for $n \geq 4$ by \cite[Theorem 5.2]{Ngu} and \cite[Lemma 6.1]{BGR2}, but was unknown for $\H^3$. As application of \eqref{HR2=hyper}--\eqref{HRhyper1}, we provide new Rellich type estimates in $\H^3$ and $\H^4$. As far as we are aware, all Rellich type inequalities with $\|\Delta_\H u\|$ in the literature were established only for $n \ge 5$.
\begin{thm}\label{prop9}
Let $n\geq 3$ and $u\in C_c^2(\H^n\backslash\{0\})$, there holds
\begin{align}
\label{HRhyper2}
\begin{split}
\|\Delta_\H u\|^2 &\geq \frac{(n-1)^4}{16}\|u\|^2 + \frac{(n-1)^2}{8}\Big\|\frac{u}{\rho}\Big\|^2 + \frac{9}{16}\Big\|\frac{u}{\rho^2}\Big\|^2\\
& \quad + \frac{(n^2-1)(n-3)(n-5)}{16}\Big\|\frac{u}{({\rm sh}\rho)^2}\Big\|^2\\
& \quad +  \frac{(n^2-1)(n-3)^2}{8}\Big\|\frac{u}{{\rm sh}\rho}\Big\|^2 + \frac{(n-1)(n-3)}{8}\Big\|\frac{u}{\rho{\rm sh}\rho}\Big\|^2.
\end{split}
\end{align}
In particular,
\begin{align*}
\|\Delta_\H u\|^2 \geq \|u\|^2 + \frac{1}{2}\Big\|\frac{u}{\rho}\Big\|^2 + \frac{9}{16}\Big\|\frac{u}{\rho^2}\Big\|^2, \quad \forall\; u\in C_c^2(\H^3\backslash\{0\}),
\end{align*}
and
\begin{align*}
\|\Delta_\H u\|^2 \geq \frac{81}{16}\|u\|^2 + \frac{9}{8}\Big\|\frac{u}{\rho}\Big\|^2 + \frac{15}{8}\Big\|\frac{u}{{\rm sh}\rho}\Big\|^2, \quad \forall\; u\in C_c^2(\H^4\backslash\{0\}).
\end{align*}
\end{thm}
The estimate \eqref{HRhyper2} improves slightly Theorem 3.1 in \cite{BGG} for $n\geq 5$.
%For more details, see $\S$ \ref{sect-hyper}.

%which deduces that for the Laplace-Beltrami operator $\Delta_\H$ and the radial Lapalcian $\Delta_{\rho,\H}$ of $\H^n$, there holds, for $n \geq 4$,
%\begin{align}
%\label{Rhyper}
%\int_{\H^n}|\Delta_\H u|^2 d v_\H \geq \int_{\H^n}|\Delta_{\rho,\H} u|^2 d v_\H, \quad \forall\; u \in C_c^2(\H^n).
%\end{align}

\medskip
It is worthy to mention that the above method to handle Hardy-Rellich inequalities by means of identities and iterations, works also for discrete situations. See the recent works \cite{HY1,SW} where this idea is applied successfully to obtain sharp discrete Hardy-Rellich inequalities in $\mathbb{N}$.

\medskip
Our paper is organized as follows. In section 2, we present some basic properties of the differential operators $T_\alpha$, $L_j$ and $R_{\alpha, k}$. The polyharmonic Hardy-Rellich identities \eqref{rHRma}, \eqref{HRm} and the Hardy-Rellich equality with general weight \eqref{HRV} are respectively proved in section 3 and 4. Finally, some examples of application are displayed in section 5, in particular we show Theorems \ref{thm6}--\ref{prop9} in $\S$ \ref{sect-hyper}.

\section{Preliminaries}
\reset

We state here some elementary but useful properties of $T_\alpha$, $L_j$ and $R_{\alpha, k}$, defined respectively in \eqref{TaLj}, \eqref{defLj} and \eqref{Rak}. By the way, we deduce the first order weighted Hardy equalities \eqref{rHR1a} and \eqref{HR1a} in Remark \ref{newrmk1}.

\begin{lem}
\label{propTa}
Let $\a, \b, \s \in \R$ and $f \in C_c^1(\R^n\backslash\{0\})$, we have
\begin{align}
\label{T0}
\left\|T_\b(f)\right\|_\a^2 = \left\|T_\a(f)\right\|_\a^2 + \frac{(\b-\a)^2}{4}\left\|f\right\|_{\a+2}^2,
\end{align}
\begin{align}
\label{newT1}
T_\b(r^\sigma f) = r^{\sigma} T_{\b-2\sigma} (f), \quad \left\|T_\b(r^\s f)\right\|_\a^2 = \left\|T_{\b -2\s}(f)\right\|_{\a-2\s}^2,
\end{align}
and
\begin{align}
\label{newT2}
\begin{split}
\langle T_\a(\p_r f),~T_{\a+2}(f)\rangle_{\a+1} = \frac{4+\a-n}{2}\|T_{\a+2} (f)\|_{\a+2}^2.
\end{split}
\end{align}
\end{lem}
Here $\langle\cdot, \cdot\rangle_\b$ denotes the inner product associated to $\|\cdot\|_\b$, that is
$$\langle f, g\rangle_\b = \int_{\R^n} \frac{f(x)g(x)}{|x|^\b}dx.$$

\medskip\noindent
{\sl Proof}. The departure point is the following simple fact: For any $\eta \in C^1_c(0, \infty)$, $M, \g \in \R$, there holds
\begin{align}
\label{0} \int_{\R_+} \Big|\eta' + \frac{M\eta}{s}\Big|^2s^\g ds = \int_{\R_+} |\eta'|^2s^\g ds - M(\g -1 - M)\int_{\R_+} |\eta|^2s^{\g - 2} ds.
\end{align}
A direct consequence is that for all $\g \in \R$ and $M', M \in \R$
\begin{align}
\label{01}
\begin{split}
& \quad \int_{\R_+} \Big|\eta' + \frac{M'\eta}{s}\Big|^2s^\g ds - \int_{\R_+} \Big|\eta' + \frac{M\eta}{s}\Big|^2s^\g ds\\
& = (M'-M)(M + M' -\g +1)\int_{\R_+} |\eta|^2s^{\g - 2} ds.
\end{split}
\end{align}

\medskip
Let $f \in C_c^1(\R^n\backslash\{0\})$. As
\begin{align}
\label{Tl}
\left\|T_\ell(f)\right\|_\a^2 = \int_{\Sp^{n-1}}\int_{\R_+} \Big|\p_r f + \frac{n-2- \ell}{2 r}f\Big|^2 r^{n-1-\a}dr d\sigma,
\end{align}
taking $M = \frac{n-2-\b}{2}$, $M' = \frac{n-2-\a}{2}$ and $\g = n-1-\a$ in \eqref{01}, we obtain readily \eqref{T0}. We get \eqref{newT1} by the definition \eqref{TaLj}. Moreover, there hold
\begin{align}
\label{newT0}
\p_r \big(T_{\a+2} (f)\big) = T_\a(\p_r f)-\frac{1}{r} T_{\a+2} (f), \quad \langle \p_r f, g\rangle_\a = -\langle f, \p_r g\rangle_\a - (n-1-\a)\langle f, g\rangle_{\a+1}.
\end{align}
\eqref{newT2} can be derived directly. \qed

\medskip
A direct application of \eqref{T0} is
\begin{align}
\label{HR13}
\|\D_r u\|_\a^2 = \left\|T_{-n}(\p_r u)\right\|_\a^2 & = \frac{(n+\a)^2}{4}\left\|\p_r u\right\|_{\a+2}^2 + \left\|T_\a(\p_r u)\right\|_\a^2,
\end{align}
which yields another classical Hardy-Rellich inequality with best constant.
\begin{align}\label{HR21r}
\|\D_r u\|_\a^2  \geq \frac{(n+\a)^2}{4}\left\|\p_r u\right\|_{\a+2}^2, \quad \forall\; u \in C_c^1(\R^n\backslash\{0\}),\; \a \in \R.
\end{align}

We recall also some well-known properties of $L_j$, see for instance \cite{MOW}.
\begin{lem}
\label{lemLj}
Let $L_j$ be given in \eqref{defLj}, there holds $|\nabla f|^2 = |\p_r f|^2 + \sum_{1\leq j \leq n} |L_j f|^2$. Moreover, we have
\begin{align}
\label{propLj}
(i)\; L_j(r)=0; \quad (ii)\; L_j \p_r=\Big(\p_r +\frac{1}{r}\Big) L_j; \quad (iii)\; \Delta = \Delta_r +\sum_{j=1}^n L_j^2;\quad (iv)\; \sum_{j=1}^n x_jL_j = 0.
\end{align}
\end{lem}

\begin{rem}
\label{newrmk1}
Taking $\ell = \alpha$, $M = \frac{n-2-\b}{2}$ and $\g = n-1-\a$ in \eqref{0}, we get readily \eqref{rHR1a}. Denote $h = r^\frac{n-2-\alpha}{2}u$, using \eqref{rHR1a}, there holds
\begin{align*}
\|\nabla u\|_\a^2 - \left(\frac{n-2 -\alpha}{2}\right)^2 \|u \|_{\a+2}^2 & = \left\|T_\a(u)\right\|_\a^2 + \sum_{j = 1}^n \|L_j u\|_\a^2\\
& = \left\|\p_r h\right\|_{n-2}^2 + \sum_{j = 1}^n \|L_j h\|_{n-2}^2 = \|\nabla h\|_{n-2}^2,
\end{align*}
which gives the Hardy equality \eqref{HR1a}.
\end{rem}

The following are some elementary properties of $R_{\a, k}$ given in \eqref{Rak}.
\begin{lem}
\label{newlemR}
For any $k\in \N$, $\a, \b \in \R$ and $f\in C_c^{k+1}(\R^n\backslash\{0\})$, there hold
\begin{align}
\label{newR1}
R_{\a,k}(r^\b f)=r^\b R_{\a-2\b, k}(f),
\end{align}
\begin{align}\label{k_a}
\|R_{\a,k} (\p_r f)\|_\a^2=\frac{(n-2k-4-\a)^2}{4}\|R_{\a+2,k} (f)\|_{\a+2}^2 +  \|R_{\a,k+1} (f)\|_{\a}^2,
\end{align}
and
\begin{align}\label{k_b}
\big\langle R_{\a,k} (\p_r f), R_{\a+2,k}(f)\big\rangle _{\a+1}=\frac{\a+2k+4-n}{2}\|R_{\a+2,k}(f)\|_{\a+2}^2.
\end{align}
\end{lem}

\medskip
\noindent
{\sl Proof}. The equality \eqref{newR1} comes directly by iterative applications of the first formula in \eqref{newT1}. We will prove \eqref{k_a}-\eqref{k_b} with induction on $k$.

\smallskip
Using successively \eqref{rHR1a} for $T_{\a+2}(f)$, \eqref{newT0} and \eqref{newT2}, we see that
\begin{align}
\label{T1}
\begin{split}
\|R_{\a, 1}(f)\|_\a^2 & = \|T_\a \circ T_{\a+2} (f)\|_\a^2\\ & = \|\p_r(T_{\a+2}(f))\|_\a^2 - \frac{(n-2-\a)^2}{4}\|T_{\a+2} (f)\|_{\a+2}^2\\
& = \left\|T_\a(\p_r f) - r^{-1} T_{\a+2} (f)\right\|_\a^2 - \frac{(n-2-\a)^2}{4} \|T_{\a+2} (f)\|_{\a+2}^2\\
& = \|T_\a(\p_r f)\|_\a^2 + \|T_{\a+2}(f)\|_{\a+2}^2 - 2\langle T_\a(\p_r f),~T_{\a+2}(f)\rangle_{\a+1}\\
&\quad - \frac{(n-2-\a)^2}{4} \|T_{\a+2} (f)\|_{\a+2}^2\\
& = \|T_\a(\p_r f)\|_\a^2 - \frac{(n-4-\a)^2}{4} \|T_{\a+2} (f)\|_{\a+2}^2.
\end{split}
\end{align}
The above equality is just \eqref{k_a} with $k = 0$, the equality \eqref{k_b} with $k=0$ is given by \eqref{newT2}.

\smallskip
Now assume that \eqref{k_a}-\eqref{k_b} hold for $k=m \in \N$.
Applying \eqref{newT0} and \eqref{newR1},
\begin{align}
\label{newR3}
\begin{split}
R_{\a,m}(\p_r T_{\a+2m+4}f)& =R_{\a,m}(T_{\a+2m+2} (\p_r f))-R_{\a,m}\left(r^{-1} T_{\a+2m+4}f\right)\\
&= R_{\a,m+1}(\p_r f) - r^{-1}R_{\a+2,m+1}(f).
\end{split}
\end{align}
On the other hand, using \eqref{k_b} with $k = m$, we have
\begin{align*}\big\langle R_{\a, m}(\p_r T_{\a+2m+4}f), R_{\a+2, m+1}(f)\big\rangle_{\a+1} & = \big\langle R_{\a, m}(\p_r T_{\a+2m+4}f), R_{\a+2, m}(T_{\a+2m+4}f)\big\rangle_{\a+1}\\
& = \frac{\a+2m+4-n}{2}\|R_{\a+2,m+1}(f)\|_{\a+2}^2.
\end{align*}
Therefore,
\begin{align}
\label{newR2}
\begin{split}
 & \quad \big\langle R_{\a, m+1}(\p_r f), R_{\a+2, m+1}(f)\big\rangle_{\a+1}\\
& = \big\langle R_{\a, m}(\p_r T_{\a+2m+4}f), R_{\a+2, m+1}(f)\big\rangle_{\a+1} +\|R_{\a+2, m+1}(f)\|_{\a+2}^2\\
& = \frac{\a+2m+4-n}{2}\|R_{\a+2,m}(T_{\a+2m+4}f)\|_{\a+2}^2+\|R_{\a+2, m+1}(f)\|_{\a+2}^2\\
& = \frac{\a+2(m+1)+4-n}{2}\|R_{\a+2,m+1}(f)\|_{\a+2}^2,
\end{split}
\end{align}
which is \eqref{k_b} for $k = m+1$. Furthermore, applying \eqref{k_a} with $k = m$ and $T_{\a+2m+4}f$, there holds
\begin{align*}
\|R_{\a, m+2} (f)\|_\a^2& =\;\|R_{\a, m+1}(T_{\a+2m+4}f)\|_\a^2\\& =\;\|R_{\a, m} (\p_r (T_{\a+2m+4}f))\|_\a^2-\frac{(n-2m-4-\a)^2}{4} \|R_{\a+2, m}(T_{\a+2m+4}f)\|_{\a+2}^2.
\end{align*}
Combining with \eqref{newR3} and \eqref{newR2}, we get
\begin{align*}
\|R_{\a, m+2} (f)\|_\a^2& =\left\|R_{\a,m+1}(\p_r f)- r^{-1}R_{\a+2,m+1}(f)\right\|_\a^2-\frac{(n-2m-4-\a)^2}{4} \|R_{\a+2, m+1} (f)\|_{\a+2}^2\\
& =\|R_{\a,m+1}(\p_r f)\|_\a^2 + \|R_{\a+2,m+1}(f)\|_{\a+2}^2\\
& \quad -2\big\langle R_{\a, m+1}(\p_r f), R_{\a+2, m+1}(f)\big\rangle_{\a+1}-\frac{(n-2m-4-\a)^2}{4} \|R_{\a+2, m+1} (f)\|_{\a+2}^2\\
& =\|R_{\a,m+1}(\p_r f)\|_\a^2 -\frac{(n-2m-6-\a)^2}{4} \|R_{\a+2, m+1}(f)\|_{\a+2}^2.
\end{align*}
This is just \eqref{k_a} with $k = m+1$, the proof is over. \qed

\section{Hardy-Rellich identities for polyharmonic operators in $\R^n$}
\reset
Here we show Theorems \ref{thmrHRma} and \ref{newthm1}. A key argument is the iterative identity for $\|R_{\a,k}(\Delta_r f)\|_\a^2$.
\begin{lem}
\label{newlemR1}
Let $\alpha \in \R$, $k \in \N$ and $f \in C_c^{k+2}(\R^n\backslash \{0\})$, then
\begin{align}
\label{iteRak}
\|R_{\a,k}(\Delta_r f)\|_\a^2 = A_{\a+ 2k+2}^2\|R_{\a+4, k} f \|_{\a+4}^2 + D_{\a+ 2k+2}\|R_{\a+2, k+1} f \|_{\a+2}^2+ \|R_{\a, k+2} f \|_{\a}^2,
\end{align}
where the coefficients $A_j$ and $D_\ell$ are given in \eqref{AD}.
\end{lem}
\noindent
{\sl Proof}. By \eqref{newR1}, we have
\begin{align*}
R_{\a,k}(\Delta_r f)=R_{\a,k}\left(r^{1-n}\p_r (r^{n-1}\p_r f)\right)=r^{1-n} R_{\a+2n-2, k}\left(\p_r (r^{n-1}\p_r f)\right).
\end{align*}
Thanks to \eqref{newR1} and \eqref{k_a}, we have
\begin{align*}
& \quad \|R_{\a,k}(\Delta_r f)\|_\a^2\\ &= \|R_{\a,k}\left(r^{1-n}\p_r (r^{n-1}\p_r f)\right)\|_\a^2\\
& = \|R_{\a+2n-2, k}\left(\p_r (r^{n-1}\p_r f)\right)\|_{\a+2n-2}^2\\
& = \frac{(n+2k+\a+2)^2}{4}\|R_{\a+2n, k}\left(r^{n-1}\p_r f\right)\|_{\a+2n}^2+ \|R_{\a+2n-2, k+1}\left(r^{n-1}\p_r f\right)\|_{\a+2n-2}^2\\
& = \frac{(n+2k+\a+2)^2}{4}\|R_{\a+2, k} (\p_r f) \|_{\a+2}^2+ \|R_{\a, k+1} (\p_r f) \|_{\a}^2\\
 &= \frac{(n+2k+\a+2)^2}{4}\left[\frac{(n-2k-\a-6)^2}{4}\|R_{\a+4, k} f \|_{\a+4}^2+ \|R_{\a+2, k+1} f \|_{\a+2}^2\right]\\
&\quad +\frac{(n-2k-\a-6)^2}{4}\|R_{\a+2, k+1} f\|_{\a+2}^2+ \|R_{\a, k+2}f \|_{\a}^2\\
& = A_{\a+ 2k+2}^2\|R_{\a+4, k} f\|_{\a+4}^2 + D_{\a+ 2k+2}\|R_{\a+2, k+1} f \|_{\a+2}^2+ \|R_{\a, k+2} f\|_{\a}^2.
\end{align*}
The proof is completed. \qed

\subsection{Proof of Theorem \ref{thmrHRma}} We will proceed by induction on $m$. For $m = 1$, by \eqref{HR13} and \eqref{T1},
\begin{align}
\label{rHRma1}
\begin{split}
\|\D_r u\|_\a^2 & = \|T_{-n}(\p_r u)\|_\a^2\\ & = \|T_\a(\p_r u)\|_\a^2 + \frac{(n+\a)^2}{4}\|\p_r u\|_{\a + 2}^2\\
& = \|R_{\a, 1}u\|_\a^2 + \frac{(n-4 - a)^2}{4}\|T_{\a+2}(u)\|_{\a + 2}^2 + \frac{(n+\a)^2}{4}\|\p_r u\|_{\a + 2}^2\\
& = A_\a^2\|u\|_{\a + 4}^2 + D_\a\|T_{\a+2}(u)\|_{\a + 2}^2  + \|R_{\a, 1}u\|_\a^2.
\end{split}
\end{align}
This is just \eqref{rHR2a}, or the identity \eqref{rHRma} with $m = 1$, $C_{0, 1, \alpha} = 1$ and $C_{1,1,\a} = D_\a > 0$.

\smallskip
Suppose that \eqref{rHRma} holds for $m \geq 1$, then
\begin{align*}
\|\D_r^{m+1} u\|_\a^2 & = \|\D_r^m(\D_r u)\|_\a^2\\ & =  \prod_{\ell = 0}^{m-1} A_{\a+4\ell}^2\times \|\D_r u\|_{\a+4m}^2+\sum_{j =0}^{2m-1} C_{j,m,\a} \left\|R_{\a+2j,2m-1-j}(\D_r u)\right\|_{\a+2j}^2.
\end{align*}
Using \eqref{rHRma1} with $\a+ 4m$,
\begin{align*}
\|\D_r u\|_{\a+4m}^2 = A_{\a+ 4m}^2\|u\|_{\a + 4m+4}^2 + D_{\a+4m}\|R_{\a+4m+2, 0} u\|_{\a + 4m+2}^2  + \|R_{\a+4m, 1}u\|_{\a+4m}^2.
\end{align*}
Moreover, applying Lemma \ref{newlemR1}, there holds, for all $0 \leq j \leq 2m -1$,
\begin{align*}
\left\|R_{\a+2j,2m-1-j}(\D_r u)\right\|_{\a+2j}^2 & = A_{\a+ 4m}^2\|R_{\a+2j + 4,2m-1-j} u\|_{\a+2j+4}^2\\
&\quad  + D_{\a+ 4m}\|R_{\a+2j+2, 2m-j} u\|_{\a+2j+2}^2+ \|R_{\a+2j,2m+1-j} u\|_{\a + 2j}^2.
\end{align*}
We conclude \eqref{rHRma} with $(m+1)$ by the iterative relation \eqref{iteC} for coefficients $C_{j, m, \a}$.
\qed

\medskip
By iterations with Lemma \ref{newlemR1}, we can claim also
\begin{prop}
\label{newcor1}
Let $n, m \geq 1$, $\a \in \R$ and $u \in C_c^{2m}(\R^n\backslash \{0\})$, there exist positive coefficients $\widetilde C_{j,m,\a}$ $0 \le j \le 2m$ such that
\begin{align}
\label{rHRmad}
\begin{split}
\|\p_r (\D_r^m u)\|_\a^2 &= \left(\frac{n-2 -\alpha}{2}\right)^2 \prod_{\ell = 0}^{m-1} A_{\a+2+4\ell}^2\times \|u\|_{\a+2+4m}^2\\
&\quad + \sum_{j =0}^{2m} \widetilde C_{j,m,\a}\left\|R_{\a+2j,2m-j}(u)\right\|_{\a+2j}^2.
\end{split}
\end{align}
\end{prop}

\medskip
\noindent{\sl Proof}. Applying \eqref{rHR1a}, there holds
\begin{align*}
\|\p_r (\D_r^m u)\|_\a^2 = \left(\frac{n-2 -\alpha}{2}\right)^2\left\|\D_r^m u\right\|_{\a+2}^2 + \left\|T_\a(\D_r^m u)\right\|_\a^2
\end{align*}
We develop the first term by \eqref{rHRma}, and the second one with the following identity.
\begin{align}
\label{rTma}
\|R_{\a, k}(\D_r^m u)\|_\a^2 & = \sum_{j = 0}^{2m} \widehat{C}_{j, m, \a, k}\|R_{\a + 2j, k+2m -j}(u)\|_{\a+2j}^2,
\end{align}
that is given by iterating \eqref{iteRak}. Here $\widehat{C}_{j, m, \a, k}$ are positive constants defined by the iterative relation as follows.
\begin{align*}
\widehat{C}_{j, m + 1, \a, k} = A_{\a+2k + 2}^2\widehat{C}_{j-2, m, \a+4, k} + D_{\a + 2k + 2}\widehat{C}_{j-1, m, \a+2, k+1} + \widehat{C}_{j, m, \a, k+2}, \quad \forall\; 0 \leq j \leq 2m+2
\end{align*}
with the convention $\widehat{C}_{0, 0, \a, k} = 1$ and $\widehat{C}_{j, m, \a, k} = 0$ if $j > 2m$ or $j < 0$. Let $\widetilde C_{j,m,\a} = \widehat C_{j,m,\a, 0}+ C_{j-1,m,\a}$, we can conclude \eqref{rHRmad}. \qed
%We can see that $$\widetilde{C}_{0, m, \a, k} = 1, \quad \widetilde{C}_{2m, m, \a, k} = \prod_{0\leq \ell \leq m-1} A_{\a + 2k + 2 + 4\ell}^2, \quad \forall \; m \geq 1.$$

\subsection{Proof of Theorem \ref{newthm1}} We begin with the second order case. Inspired by \cite{MOW} (with $\a = 0$), we prove first \eqref{HR2a}. In fact,
\begin{align}
\label{newHR2a}
\|\D u\|_\a^2 & =\Big\|\D_r u + \sum_{j=1}^n L_j^2 u\Big\|_\a^2 =\; \|\D_r u\|_\a^2 + \Big\|\sum_{j=1}^n L_j^2 u\Big\|_\a^2 + 2\Big\langle \D_r u, \sum_{j=1}^n L_j^2 u\Big\rangle_\a.
\end{align}
Applying formula $(iv)$ in \eqref{propLj}, direct calculation yields
\begin{align}
\label{ibpa}
\Big\langle h, \sum_{j=1}^n L_j^2 u\Big\rangle_\a= - \sum_{j=1}^n \langle L_j h, L_j u\rangle_\a + (n-1) \Big\langle h, \sum_{j=1}^n x_jL_j u\Big\rangle_{\a+2} = - \sum_{j=1}^n \langle L_j h, L_j u\rangle_\a.
\end{align}
Moreover, by \eqref{propLj},
\begin{align*}
L_j(\Delta_r u)= L_j \Big(r^{1-n} \p_r(r^{n-1}\p_r u)\Big) & = r^{1-n} \Big(\p_r + \frac{1}{r}\Big) L_j(r^{n-1}\p_r u)\\
& = r^{1-n} \Big(\p_r + \frac{1}{r}\Big)\left[r^{n-1} \Big(\p_r + \frac{1}{r}\Big) L_j u\right]\\
& = \p_r^2 (L_j u) +\frac{n+1}{r} \p_r (L_j u) +\frac{n-1}{r^2} L_j u.
\end{align*}
Combining with \eqref{newT0}, \eqref{ibpa} and \eqref{rHR1a}, we have then
\begin{align*} \Big\langle \D_r u, \sum_{j=1}^n L_j^2 u\Big\rangle_\a
& = -\sum_{j=1}^n \langle L_j(\D_r u), L_j u\rangle_\a\\
&= - \sum_{j=1}^n \Big\langle \p_r^2 (L_j u) +\frac{n+1}{r} \p_r (L_j u) +\frac{n-1}{r^2} L_j u, L_j u\Big\rangle_\a\\
& = -\sum_{j=1}^n\Big[ (n-1)\| L_j u\|_{\a+2}^2-\|\p_r(L_j u)\|_\a^2 + (2+\alpha)\langle \p_r(L_j u), L_j u\rangle_{\a+1}\Big]\\
& = A_\a \sum_{j=1}^n \| L_j u\|_{\a+2}^2 + \sum_{j=1}^n \|T_\a (L_j u)\|_\a^2.
\end{align*}
Inserting this into \eqref{newHR2a}, we obtain \eqref{HR2a}, hence \eqref{HRm} holds true for $m = 1$ seeing \eqref{rHR2a} or \eqref{rHRma}.

\medskip
Now the Hardy-Rellich identity \eqref{HRm} can be derived by induction. Suppose that \eqref{HRm} holds true for $\Delta^k$, $1 \le k \le m$. Applying successively \eqref{HR2a} and \eqref{rHR2a},  or equivalently \eqref{HRm} with $m = 1$,
there holds
\begin{align*}
\|\D^{m+1} u\|_\a^2
%& =\|\D_r(\D^m u)\|_\a^2+ \left\|(\Delta - \D_r) (\D^m u)\right\|_\a^2\\&\quad + 2A_\a\sum_{j=1}^n\left\|L_j(\D^m u)\right\|_{\a+2}^2 + 2\sum_{j=1}^n\left\|T_\a\left(L_j (\D^m u)\right)\right\|_\a^2\\
& = A_\a^2\left\|\D^m u\right\|_{\a+4}^2 + D_\a\left\|T_{\a+2}(\D^m u)\right\|_{\a+2}^2 + \left\|R_{\a, 1}(\D^m u)\right\|_\a^2\\
& \quad + \left\|(\Delta - \D_r) (\D^m u)\right\|_\a^2 + 2A_\a\sum_{j=1}^n\left\|L_j(\D^m u)\right\|_{\a+2}^2 + 2\sum_{j=1}^n\left\|T_\a\left(L_j (\D^m u)\right)\right\|_\a^2.
\end{align*}
The above equality yields readily \eqref{HRm} for $\Delta^{m+1}$ with the identity for $\|\Delta^m u\|_{\a+4}$ and the definition of coefficients $H_{k, \alpha}$.
%\footnote{\color{red} There were some repeats for the statement of Theorem, it's better to write differently}
\qed

\begin{rem}
\label{remThHRm}
Notice that the coefficients $A_\beta$ are not always nonnegative. Considering the sign of $A_\b$, we see the following well known Hardy-Rellich inequality: Let $n > 2m$, $-n < \a < n - 4m$, there holds
\begin{align}
\label{ineHR2ma}
\|\D^m u\|_\a^2 \geq \prod_{\ell = 0}^{m-1} A_{\a+4\ell}^2 \times \|u\|_{\a+4m}^2, \quad \forall\; u \in C_c^{2m}(\R^n\backslash\{0\}).
\end{align}
Consequently, for $n > 2m$, $-n < \a + 2 < n - 4m$, we have
\begin{align}
\label{ineHR(2m+1)}
\|\nabla \D^m u\|_\a^2 \geq \|\p_r \D^m u\|_\a^2 \geq \frac{(n - 2 - \a)^2}{4} \prod_{\ell = 0}^{m-1} A_{\a+4\ell + 2}^2 \times \|u\|_{\a+4m+2}^2, \quad \forall\; u \in C_c^{2m+2}(\R^n\backslash\{0\}).
\end{align}
Clearly, positive remainder terms can be added to the above inequalities by equalities \eqref{HRm} and \eqref{HR1a}. We may also use \eqref{HRV} to get diverse formulae or relax the assumption on $\a$ using \eqref{Sphere} below.
\end{rem}

\section{Higher order Hardy-Rellich identities with general weights}
\reset
This section is devoted to Theorem \ref{thmHRV}. To be more transparent, we consider first the even order cases. When $k = 2m$, let us rewrite \eqref{HRV} as follows.
\begin{align}
\label{HRV2m}
\begin{split}
&\quad \int_{\M} V|\Delta_g^m u|^2 dg-\int_{\M}\frac{\Delta_g^m (V\Delta_g^m f)}{f}u^2 dg\\& =-2\sum_{\ell=0}^{m-1}\int_{\M}\frac{\Delta_g^{m-1-\ell}(V\Delta_g^m f)}{\Delta_g^\ell f}\Big\|\nabla_g(\Delta_g^\ell u)-\frac{\Delta_g^\ell u}{\Delta_g^\ell f} \nabla_g(\Delta_g^\ell f)\Big\|^2 dg \\
&\quad +\sum_{\ell=1}^m \int_{\M}\frac{\Delta_g^{m-\ell}(V\Delta_g^m f)}{\Delta_g^\ell f}\Big|\Delta_g^\ell u -\frac{\Delta_g^{\ell-1} u}{\Delta_g^{\ell-1} f}\Delta_g^\ell f \Big|^2 dg.
\end{split}
\end{align}
Observe that the last term of \eqref{HRV} is inserted in the previous line with $\ell = m$.

\medskip
We use the induction to show \eqref{HRV2m}. For any $\phi\in C^2(\M)$, $u \in C_c^2(\M)$, there holds
\begin{align*}
\int_{\M} V|\Delta_g u +  \phi u|^2 dg
& =\int_{\M} V\Big[|\Delta_g u|^2 + \phi^2 u^2 + 2 \phi u\Delta_g u\Big]dg\\
& =\int_{\M} V|\Delta_g u|^2dg +\int_{\M} V\phi^2 u^2dg -2\int_{\M} V\phi|\nabla_g u|^2dg + \int_{\M} \Delta_g(V \phi)u^2dg\\
& =\int_{\M} V|\Delta_g u|^2dg +\int_{\M} V\phi^2 u^2dg +2\int_{\M} \frac{{\rm div_g}(V\phi\nabla_g f)}{f} u^2dg\\
& \quad + \int_{\M} \Delta_g(V \phi)u^2 dg - 2 \int_{\M} V\phi\Big\|\nabla_g u - \frac{u}{f}\nabla_g f\Big\|^2 dg.
\end{align*}
This implies that
\begin{align*}
\int_{\M} V|\Delta_g u|^2dg& =\int_{\M}\Big[-2 \frac{{\rm div_g}(V\phi\nabla_g f)}{f}-V\phi^2- \Delta_g(V \phi)\Big] u^2 dg\\
&\quad +2 \int_{\M} V\phi\Big\|\nabla_g u - \frac{u}{f}\nabla_g f\Big\|^2 dg +\int_{\M} V|\Delta_g u +  \phi u|^2 dg.
\end{align*}
Set $\phi:=-\frac{\Delta_g f}{f}$ with $f\in C^2$, then
\begin{align*}
2 \frac{{\rm div_g}(V\phi\nabla_g f)}{f}+V\phi^2+ \Delta_g(V \phi) & = 2\frac{\nabla_g(V\phi)\cdot\nabla_g f}{f} + V\phi \frac{\Delta_g f}{f} +\Delta_g (V\phi)\\
 & =\frac{\Delta_g(V\phi f)}{f}\\
& = -\frac{\Delta_g (V\Delta_g f)}{f}.
\end{align*}
Consequently
\begin{align}
\label{Delta1}
\begin{split}
\int_{\M} V|\Delta_g u|^2dg & =\int_{\M} \frac{\Delta_g (V\Delta_g f)}{f} u^2 dg -2 \int_{\M} V\frac{\Delta_g f}{f}\Big\|\nabla_g u - \frac{u}{f}\nabla_g f\Big\|^2 dg\\
& \quad +\int_{\M} V\Big|\Delta_g u -\frac{\Delta_g f}{f} u\Big|^2 dg,
\end{split}
\end{align}
i.e. we get \eqref{HRV2m} for $m=1$. Assume that \eqref{HRV2m} holds for some $m \in \N$, let us consider the order $(m+1)$.
Using \eqref{Delta1} with $f = \psi$ and $\Delta_g^mu$,
\begin{align}
\label{new4.3}
\begin{split}
&\quad \int_{\M} V|\Delta_g^{m+1} u|^2 dg - \int_{\M}\frac{\Delta_g(V\Delta_g\psi)}{\psi} |\Delta_g^m u|^2 dg\\
& = - 2 \int_{\M} V\frac{\Delta_g\psi}{\psi} \Big\|\nabla_g(\Delta_g^m u) - \Delta_g^m u \frac{\nabla_g \psi}{\psi}\Big\|^2 dg + \int_{\M} V\Big|\Delta_g^{m+1}u - \frac{\Delta_g \psi}{\psi} \Delta_g^m u\Big|^2dg.
\end{split}
\end{align}
Moreover, applying \eqref{HRV2m} with $m$ and $V_m = \frac{\Delta_g(V\Delta_g\psi)}{\psi}$ for the second integral in \eqref{new4.3}, we obtain
\begin{align*}
&\quad \int_{\M} V|\Delta_g^{m+1} u|^2 dg - \int_{\M}\frac{\Delta_g^m(V_m\Delta_g^m f)}{f} u^2 dg\\
& = -2\sum_{\ell=0}^{m-1}\int_{\M}\frac{\Delta_g^{m-1-\ell}(V_m\Delta_g^m f)}{\Delta_g^\ell f}\Big\|\nabla_g(\Delta_g^\ell u)-\frac{\Delta_g^\ell u}{\Delta_g^\ell f} \nabla_g(\Delta_g^\ell f)\Big\|^2 dg\\
&\quad +\sum_{\ell=1}^m \int_{\M}\frac{\Delta_g^{m-\ell}(V_m\Delta_g^m f)}{\Delta_g^\ell f}\Big|\Delta_g^\ell u -\frac{\Delta_g^{\ell-1} u}{\Delta_g^{\ell-1} f}\Delta_g^\ell f \Big|^2 dg\\
&\quad - 2 \int_{\M} V\frac{\Delta_g\psi}{\psi} \Big\|\nabla_g(\Delta_g^m u) - \Delta_g^m u \frac{\nabla_g \psi}{\psi}\Big\|^2 dg + \int_{\M} V\Big|\Delta_g^{m+1}u - \frac{\Delta_g \psi}{\psi} \Delta_g^m u\Big|^2dg.
\end{align*}
Let now $\psi = \Delta_g^m f$ so that $V_m\Delta_g^m f = \Delta_g(V\Delta_g^{m+1}f)$ and $\frac{\Delta_g\psi}{\psi} = \frac{\Delta_g^{m+1}f}{\Delta_g^m f}$, we arrive at
\begin{align*}
&\quad \int_{\M} V|\Delta_g^{m+1} u|^2 dg - \int_{\M}\frac{\Delta_g^{m+1} (V\Delta_g^{m+1}f)}{f} u^2 dg\\
& = -2\sum_{\ell=0}^m\int_{\M}\frac{\Delta_g^{m-\ell}(V\Delta_g^{m+1} f)}{\Delta_g^\ell f}\Big\|\nabla_g(\Delta_g^\ell u)-\frac{\Delta_g^\ell u}{\Delta_g^\ell f} \nabla_g(\Delta_g^\ell f)\Big\|^2 dg \\
&\quad +\sum_{\ell=1}^{m+1}\int_{\M}\frac{\Delta_g^{m+1-\ell}(V\Delta_g^{m+1} f)}{\Delta_g^\ell f}\Big|\Delta_g^\ell u -\frac{\Delta_g^{\ell-1} u}{\Delta_g^{\ell-1} f}\Delta_g^\ell f \Big|^2 dg.
\end{align*}
This means that \eqref{HRV2m} holds true for any $m \geq 1$.

\medskip
Next we consider odd order cases of \eqref{HRV}, i.e.~$k = 2m+1$. When $k=1$, \eqref{HRV} is just the identity \eqref{HYequa1}. Consider $k=2m+1$ for general $m \in \N$. Using \eqref{HYequa1} with $\Delta_g^m u$ and $f = \psi$, there holds
\begin{align*}
\int_{\M} V\|\nabla_g\big(\Delta_g^m u\big)\|^2 dg = -\int_{\M} \frac{{\rm div_g}(V\nabla_g\psi)}{\psi} |\Delta_g^m u|^2 dg +\int_{\M} V\Big\|\nabla_g\big(\Delta_g^m u\big) - \frac{\Delta_g^m u}{\psi}\nabla_g \psi\Big\|^2dg.
\end{align*}
Applying \eqref{HRV2m} with the weight $V_m :=  -\frac{{\rm div_g}(V\nabla_g\psi)}{\psi}$, there holds
\begin{align*}
&\quad \int_{\M} V\|\dD_{2m+1} u\|^2 dg - \int_{\M}\frac{\Delta_g^m (V_m\Delta_g^m f)}{f} u^2 dg\\
& = -2\sum_{\ell=0}^{m-1}\int_{\M}\frac{\Delta_g^{m-1-\ell}(V_m\Delta_g^m f)}{\Delta_g^\ell f}\Big\|\nabla_g(\Delta_g^\ell u)-\frac{\Delta_g^\ell u}{\Delta_g^\ell f} \nabla_g(\Delta_g^\ell f)\Big\|^2 dg \\
&\quad +\sum_{\ell=1}^m\int_{\M}\frac{\Delta_g^{m-\ell}(V_m\Delta_g^m f)}{\Delta_g^\ell f}\Big|\Delta_g^\ell u -\frac{\Delta_g^{\ell-1} u}{\Delta_g^{\ell-1} f}\Delta_g^\ell f \Big|^2 dg +\int_{\M} V\Big\|\nabla_g\big(\Delta_g^m u\big) - \frac{\Delta_g^m u}{\psi}\nabla_g \psi\Big\|^2dg.
\end{align*}
Choose $\psi = \Delta_g^m f = \dD_{2m}f$ so that $\Delta_g^j(V_m\Delta_g^m f) = \dD^*_{2j+1}(V\dD_{2m+1} f)$ and $\nabla_g \psi = \dD_{2m+1} f$. We deduce then
\begin{align*}
&\quad \int_{\M} V\|\dD_{2m+1} u\|^2 dg - \int_{\M}\frac{\dD^*_{2m+1}(V\dD_{2m+1} f)}{f} u^2 dg\\
 & = -2\sum_{\ell=0}^{m-1}\int_{\M}\frac{\dD^*_{2m-1-2\ell}(V\dD_{2m+1} f)}{\dD_{2\ell} f}\Big\|\dD_{2\ell+1} u -\frac{\dD_{2\ell}u}{\dD_{2\ell} f}\dD_{2\ell+1}f \Big\|^2 dg \\
&\quad  + \sum_{\ell=1}^m \int_{\M}\frac{\dD^*_{2m+1-2\ell}(V\dD_{2m+1} f)}{\dD_{2\ell} f}\Big|\dD_{2\ell} u -\frac{\dD_{2\ell}u}{\dD_{2\ell-2} f}\dD_{2\ell}f\Big|^2 dg\\
&\quad + \int_{\M} V\Big\|\dD_{2m+1} u - \frac{\dD_{2m} u}{\dD_{2m} f}\dD_{2m+1}f\Big\|^2dg.
\end{align*}
This is just the formula \eqref{HRV} with $k = 2m+1$. So we are done. \qed

\medskip
Remark that $\dD_{k - 2m}^*\circ\dD_k = (-1)^k\Delta_g^{k-m}$ for any $k, m \in \N$ and $2m < k$. Let $V \equiv 1$ in \eqref{HRV}, we have
\begin{align}
\label{HRV=1}
\begin{split}
&\quad \int_{\M} \|\dD_k u\|^2 dg-\int_{\M}\frac{(-\Delta_g)^k f}{f}u^2 dg\\
& = 2\sum_{\ell \ge 0, 2\ell\leq k-2}\int_{\M}\frac{(-\Delta_g)^{k-1-\ell} f}{(-\Delta_g)^\ell f}\Big\|\dD_{2\ell+1} u -\frac{\dD_{2\ell}u}{\dD_{2\ell} f}\dD_{2\ell+1}f \Big\|^2 dg \\ & \quad +\sum_{\ell \geq 1, 2\ell \le k-1}\int_{\M}\frac{(-\Delta_g)^{k-\ell} f}{(-\Delta_g)^\ell f}
\Big|\dD_{2\ell} u -\frac{\dD_{2\ell}u}{\dD_{2\ell-2} f}\dD_{2\ell}f\Big|^2 dg\\
& \quad + \int_\M \Big\|\dD_k u - \frac{\dD_{2[(k-1)/2]}u}{\dD_{2[(k-1)/2]}f}\dD_k f\Big\|^2 dg.
\end{split}
\end{align}
Corollary \ref{coriHRV} is easily deduced by \eqref{HRV=1}. We can also apply \eqref{HRV} with $V = |x|^{-\a}$ and $f = |x|^{\frac{\a+2m-n}{2}}$ to refind \eqref{ineHR2ma} or \eqref{ineHR(2m+1)}, we leave the detail for interested readers.

\section{Some applications}
\reset
Here we show several applications where we apply identities to prove some old and new Hardy-Rellich inequalities. As already mentioned, we focus on the second order cases which provide already convincing examples.

\subsection{Between Laplacian and radial Laplacian}
For general $\a \in \R$ and $m \ge 1$, we cannot always compare easily the weighted norms $\|\D^mu\|_\a$ and $\|\D_r^mu\|_\a$. Take $n \geq 2$, $m = 1$ and $u \in C_c^2(\R^n\backslash\{0\})$, by \eqref{HR2a},
\begin{align}
\label{LrL}
\|\D u\|_\a^2 - \|\D_r u\|_\a^2 & = 2A_\a\sum_{j=1}^n\left\|L_ju\right\|_{\a+2}^2 + \left\|(\Delta - \D_r) u\right\|_\a^2 + 2\sum_{j=1}^n\left\|T_\a\left(L_ju\right)\right\|_\a^2.
\end{align}
When $A_\alpha$ is nonnegative, i.e.~$-n \le \a \le n-4$, we have obviously $\|\D u\|_\a^2 \geq \|\D_r u\|_\a^2$. The range of $\a$ to ensure $\|\D u\|_\a \geq \|\D_r u\|_\a$ can be extended by the following simple fact.
\begin{lem}
\label{lemSphere}
For any $\a \in \R$, $f\in C_c^2(\R^n \backslash \{0\})$,
\begin{align}
\label{Sphere}
\left\|(\Delta - \D_r) f\right\|_\a^2 \geq (n-1)\sum_{j=1}^n\|L_j f\|_{\a+2}^2.
\end{align}
\end{lem}

\medskip
\noindent
{\sl Proof.} Let $-\D_\sigma$ be the spherical Laplacian over the unit sphere $\Sp^{n-1} \subset \R^n$. Let $h(\sigma) = f(r\sigma)$ be defined on $\Sp^{n-1}$ with $r > 0$, then
\begin{align*}
\int_{\Sp^{n-1}} |(\Delta - \D_r)f|^2(r\sigma) d\sigma & = \frac{1}{r^4}\int_{\Sp^{n-1}} |\Delta_\sigma h|^2 d\sigma\\
& \geq \frac{n- 1}{r^4}\int_{\Sp^{n-1}}|\nabla_\sigma h|^2 d\sigma = \frac{n- 1}{r^2}\int_{\Sp^{n-1}}\sum_{j= 1}^n |L_j f|^2(r\sigma) d\sigma,
\end{align*}
because $(n-1)$ is the first eigenvalue of $-\D_\sigma$ in $\dot{H}^1(\Sp^{n-1})$, the subspace of functions in $H^1(\Sp^{n-1})$ having null average. \qed

\smallskip
Using \eqref{LrL} and \eqref{Sphere}, we see that
\begin{align*}
%\label{LrL}
\|\D u\|_\a^2 - \|\D_r u\|_\a^2 \geq (2A_\a + n-1)\sum_{j=1}^n\left\|L_ju\right\|_{\a+2}^2.
\end{align*}
Hence $\|\D u\|_\a^2 \geq \|\D_r u\|_\a^2$ if $|\alpha+2| \leq \sqrt{n^2 -2n + 2}$.

\subsection{Around some optimal Rellich inequalities}
\label{sect5.2}
For similar reason as above, the Rellich inequality
\begin{align}
\label{HR21bis}
\|\D u\|_\a^2 \geq \frac{(n+\a)^2}{4} \left\|\nabla u\right\|^2_{\a+2}
\end{align}
does not hold in $C_c^2(\R^n\backslash\{0\})$ for all $\a \in \R$, contrary to the radial derivative case (see \eqref{HR21r}). As already mentioned, Beckner and Tertikas-Zographopoulos \cite{B, TZ} showed that \eqref{HR21bis} holds true for all $u \in C_c^2(\R^n\backslash\{0\})$
if and only if $$|3\a + n+4| \le 2\sqrt{n^2 - n+1},$$
More precisely, the best constants
$$E_{n, \a} = \inf_{u\in C_c^2(\R^n\backslash\{0\}), u \ne 0}\frac{\|\D u\|_\a^2}{\|\nabla u\|_{\a+2}^2}$$ are known for general $\alpha$ and $n$, see \cite{B,TZ}. In particular, $E_{n,0} = \frac{n^2}{4}$ for $n\ge 5$, $E_{4, 0} = 3$ and $E_{3, 0} = \frac{25}{36}$, $E_{2, 0} = 0$. Here we just explain how to find quickly the optimal lower bound for $E_{n, 0}$ using Hardy-Rellich identities.

\smallskip
Indeed, combining equalities \eqref{HR2a}, \eqref{HR13} and estimate \eqref{Sphere} with $\a = 0$,
\begin{align}
\label{HR21a}
\begin{split}
&\quad \|\D u\|^2 - \frac{n^2}{4} \|\nabla u\|^2_2\\  & = \|T_0(\p_r u)\|^2 + \Big[\frac{n(n-4)}{2}-\frac{n^2}{4}\Big]\sum_{j=1}^n \left\|L_j u\right\|^2_2 + \left\|(\Delta - \D_r) u\right\|^2+2\sum_{j=1}^n\|T_0(L_j u)\|^2\\
& \geq \frac{n^2 - 4n - 4}{4}\sum_{j=1}^n \left\|L_j u\right\|^2_2.
\end{split}
\end{align}
Hence $E_{n,0} \ge \frac{n^2}{4}$ if $n\ge 5$ and $E_{4, 0} \ge 3$. To handle the $\R^3$ case, we decompose $u = \ou + u_\sigma$ where $\overline u$ is the spherical average function of $u$. So $\ou$ is radial symmetric,
\begin{equation*}
\|\D \ou\|^2 = \|\D_r \ou\|^2 \geq \frac{9}{4}\|\p_r \ou\|_2^2 = \frac{9}{4}\|\nabla \ou\|_2^2, \quad \forall\; u \in C_c^2(\R^3\backslash\{0\}).
\end{equation*}
On the other hand, using  \eqref{LrL} and \eqref{HR21r} with $u_\s$, $n=3$ and $\a = 0$, there holds
\begin{align}
\label{HR3T}
\begin{split}
\|\D u_\sigma\|^2 - C\|\nabla u_\sigma\|_2^2 & = \|\D_r u_\sigma\|^2 - C\|\p_r u_\sigma\|_2^2 -\Big(\frac{3}{2}+ C\Big)\sum_{j=1}^3\left\|L_ju_\sigma\right\|_2^2\\
&\quad  + \left\|(\Delta - \D_r) u_\sigma\right\|^2 + 2\sum_{j=1}^3\left\|T_0(L_ju_\sigma)\right\|^2\\
& \ge \Big(1 - \frac{4C}{9}\Big)\|\D_r u_\sigma\|^2 -\Big(\frac{3}{2}+ C\Big)\sum_{j=1}^3\left\|L_ju_\sigma\right\|_2^2 + \left\|(\Delta - \D_r) u_\sigma\right\|^2.
\end{split}
\end{align}
Let $\{\mu_k\}$ be the increasing sequence of eigenvalues for $-\Delta_\s = \D_r - \Delta$ in $\dot{H}^1({\mathbb S}^2)$, we know that $\mu_1 = 2$ and $\mu_2 = 6$. Seeing \eqref{HR3T}, by decomposition, to have $\|\D u_\sigma\|^2 - C\|\nabla u_\sigma\|_2^2\ge 0$ with $C \leq \frac{9}{4}$, we need only to check the component with $\mu_1$. Since $u_\sigma(r\cdot), L_ju_\sigma(r\cdot) \in \dot{H}^1({\mathbb S}^2)$ for any $r > 0$, we have
$$\left\|(\Delta - \D_r) u_\sigma\right\|^2 \ge 2\|\nabla_\sigma u_\sigma\|_2^2 \ge 4\left\|u_\sigma\right\|_4^2.$$
So it is enough to check
\begin{align}
\label{HR3T3}
\Big(1 - \frac{4C}{9}\Big)\|\D_r f\|^2 - (3+2C)\left\|f\right\|_4^2 +  4\left\|f\right\|_4^2\ge 0, \quad \forall\; f \in C_c^2(\R^3\backslash\{0\}).
\end{align}
By Theorem 1.1 in $\R^3$, $\|\D_r f\|^2 \ge \frac{9}{16}\|f\|_4^2$, we can claim \eqref{HR3T3} hence
 $\|\D u_\sigma\|^2 - C\|\nabla u_\sigma\|_2^2\ge 0$ if
$$\Big(1 - \frac{4C}{9}\Big)\frac{9}{16} + 1 - 2C \geq 0, \quad \mbox{equally $C \le \frac{25}{36}$.}$$
Applying Fubini Theorem, it's not difficult to see that $\langle \D\ou, \D u_\sigma\rangle = 0$, thus we get $E_{3, 0} \geq \frac{25}{36}$ by
$$\|\D u\|^2 = \|\D \ou\|^2 +  \|\D u_\sigma\|^2 \geq \frac{9}{4}\|\nabla \ou\|_2^2+\frac{25}{36}\|\nabla u_\sigma\|_2^2 \geq \frac{25}{36}\|\nabla u\|_2^2. $$

\medskip
It is also known that the optimal Rellich inequality $\|\D u\|_\a^2 \geq A_\a^2 \|u\|^2_{\a+4}$ holds true in $C_c^2(\R^n\backslash\{0\})$ if and only if $|\a+2| \le \sqrt{n^2 - 2n +2}$, see \cite[Theorem 4.1]{CM}. Under suitable condition on $\alpha$, we can give some explicit and simple remainder terms for
$$\|\D u\|^2_\a - \frac{(n+\a)^2}{4}\|\nabla u\|_{\a+2}^2 \quad\mbox{and}\quad \|\D u\|^2_\a - A_\a^2\|u\|_{\a+4}^2,$$
always by means of identities. The following result generalizes Theorems 2.4 in \cite{TZ} with $\a = 0$.
\begin{prop}
\label{propTZ1}
Let $\a \in \R$ satisfy $|\a+2| \le \sqrt{n-1}$. Then for any $u \in C_c^2(\R^n\backslash\{0\})$,
\begin{align*}
& \quad \|\D u\|_\a^2 - \frac{(n+\a)^2}{4} \left\|\nabla u\right\|^2_{\a+2} - \frac{(n-4-\a)^2}{4}\left\|\nabla \Big(|x|^{\frac{n-4-\a}{2}}u\Big)\right\|_{n-2}^2\\
& = \|\D u\|_\a^2 - A_\a^2 \|u\|^2_{\a+4}- D_\a\left\|\nabla \Big(|x|^{\frac{n-4-\a}{2}}u\Big)\right\|_{n-2}^2 \geq 0.
\end{align*}
\end{prop}

\medskip
\noindent
{\sl Proof.} Take $\phi =|x|^{\frac{n-4-\a}{2}}u$, $\a = n-2$ in \eqref{rHR1a},
\begin{align}\label{HR27}
\|\nabla \phi\|_{n-2}^2= \|T_{n-2} (\phi)\|_{n-2}^2 +\sum_{j=1}^n \|L_j \phi \|^2_{n-2} = \|T_{\a+2} u\|_{\a+2}^2 +\sum_{j=1}^n \|L_j u\|_{\a+2}^2.
\end{align}
Combining with \eqref{HR2a}, developing $\|T_\a(\p_r u)\|_\a^2$ by \eqref{k_a} and using \eqref{Sphere}, there holds
\begin{align*}
& \quad \|\D u\|_\a^2 - \frac{(n+\a)^2}{4} \left\|\nabla u\right\|^2_{\a+2} -\frac{(n-4-\a)^2}{4}\|\nabla \phi\|_{n-2}^2\\& = -(2+\a)^2 \sum_{j=1}^n \left\|L_j u\right\|^2_{\a+2} + \left\|(\Delta - \D_r) u\right\|_\a^2+ \|R_{\a, 1}(u)\|_{\a}^2 + 2\sum_{j=1}^n\|T_\a(L_j u)\|_\a^2\\
& \geq \Big[n-1 -(2+\a)^2\Big] \sum_{j=1}^n \left\|L_j u\right\|^2_{\a+2}.
\end{align*}
On the other hand, applying Theorem \ref{newthm1} with $m = 1$ and \eqref{HR27}, we have
\begin{align*}
&\quad \|\D u\|_\a^2 - A_\a^2 \|u\|^2_{\a+4}- D_\a\|\nabla \phi\|_{n-2}^2\\& =  -(2+\a)^2 \sum_{j=1}^n \left\|L_j u\right\|^2_{\a+2} + \left\|(\Delta - \D_r) u\right\|_\a^2+ \|R_{\a, 1}(u)\|_{\a}^2 + 2\sum_{j=1}^n\|T_\a(L_j u)\|_\a^2
\end{align*}
Therefore, the left hand side of the two above formulas are identical, the proof is completed. \qed
%the following result which generalizes Theorem 2.4 in \cite{TZ} with $\a = 0$

\medskip
Consider still $\phi =|x|^{\frac{n-4-\a}{2}}u$. As $A_{n-4} = 0$ and $D_{n-4} = (n-2)^2$, applying \eqref{HRm} with $m = 1$, direct calculation yields
\begin{align*}
%\label{HR28}
%\begin{split}
\|\D \phi\|_{n-4}^2
 & =  (n-2)^2\|T_{n-2}(\phi)\|_{n-2}^2 + \|R_{n-4, 1}(\phi)\|_{n-4}^2\\
&\quad + \left\|(\Delta - \D_r) \phi\right\|_{n-4}^2 + 2\sum_{j=1}^n\|T_{n-4}(L_j\phi)\|_{n-4}^2\\
& = (n-2)^2\|T_{\a+2} (u)\|_{\a+2}^2 + \|R_{\a, 1} (u) \|_{\a}^2 + \left\|(\Delta - \D_r) u\right\|_{\a}^2 + 2\sum_{j=1}^n\left\|T_{\a}\left(L_j u \right)\right\|_{\a}^2.
%\end{split}
\end{align*}
Combining \eqref{HR21a}, \eqref{k_a} and \eqref{Sphere}, we have
\begin{align}
\label{new5.8}
\begin{split}
&\quad \|\D u\|_\a^2 - \frac{(n+\a)^2}{4} \left\|\nabla u\right\|^2_{\a+2}- \frac{(n-4-\a)^2}{4(n-2)^2}\|\Delta \phi\|_{n-4}^2\\
& = \Big[1 - \frac{(n-4-\a)^2}{4(n-2)^2}\Big]\Big(\|R_{\a, 1} (u) \|_\a^2 + \left\|(\Delta - \D_r) u\right\|_{\a}^2 + 2 \sum_{j=1}^n \|T_{\a} (L_j (u))\|_{\a}^2\Big) \\
&\quad + \Big[2A_\a - \frac{(n+\a)^2}{4}\Big]\sum_{j=1}^n\left\|L_j u\right\|_{\a+2}^2\\
& \geq \Big[1 - \frac{(n-4-\a)^2}{4(n-2)^2}\Big]\Big(\|R_{\a, 1} (u) \|_\a^2 + 2 \sum_{j=1}^n \|T_{\a} (L_j (u))\|_{\a}^2\Big)\\
& \quad + \left\{\Big[1 - \frac{(n-4-\a)^2}{4(n-2)^2}\Big](n - 1) + 2A_\a - \frac{(n+\a)^2}{4}\right\}\sum_{j=1}^n\left\|L_j u\right\|_{\a+2}^2.
\end{split}
\end{align}
Hence the left hand side of \eqref{new5.8} is nonnegative if
$$4(n-2)^2 - (n-4-\a)^2 \geq 0 \quad \mbox{and}\quad \Big[1 - \frac{(n-4-\a)^2}{4(n-2)^2}\Big](n - 1) + 2A_\a - \frac{(n+\a)^2}{4} \geq 0,$$
or equivalently
\begin{align}
\label{TZ2c}
n \geq 3,\quad -n \leq \a \leq \frac{n^3 - 9n^2 + 25n -24}{3n^2 -11n+11}.
\end{align}
So we extend Theorem 2.7 in \cite{TZ} as follows.
\begin{prop}
\label{propTZ2}
Let $n, \a$ satisfy \eqref{TZ2c}, then for any $u \in C_c^2(\R^n\backslash\{0\})$,
\begin{align}
\label{TZ2}
&\quad \|\D u\|_\a^2 - \frac{(n+\a)^2}{4} \left\|\nabla u\right\|^2_{\a+2} \geq \frac{(n-4-\a)^2}{4(n-2)^2}\left\|\D \Big(|x|^{\frac{n-4-\a}{2}} u \Big)\right\|_{n-2}^2.
\end{align}
\end{prop}
%As \eqref{TZ2} holds true for $n \geq 5$ and $\a = 0$, the above result generalizes Theorem 2.7 in \cite{TZ}, but we cannot claim that the condition \eqref{TZ2c} is optimal to have \eqref{TZ2}.

\subsection{Hardy-Rellich identities with weights involving logarithmic}
\label{sect-ln}
%\reset
Here we show that the iterative idea to getting identities \eqref{rHRma} and \eqref{HRm} can be applied for weights $|x|^\a\ln|x|$. As consequence, we will prove Theorem \ref{thmHRln}. Let $f \in C_c^1(\R^n\backslash (\{0\}\cup \Sp^{n-1}))$. Denote
$$G_{\a,\b}(f):=T_\a(f)+\frac{2\b-1}{2r|\ln r|}f, \quad \forall\; \a, \b \in \R.$$
Direct calculation implies that for any $f \in C_c^1(\R^n\backslash(\{0\}\cup\Sp^{n-1}))$,
\begin{align}\label{HR32}
\left\||\ln r|^\b T_\alpha(f)\right\|_\a^2 = \frac{(2\b-1)^2}{4}\left\||\ln r|^{\b-1}f\right\|_{\a+2}^2 + \||\ln r|^\b G_{\alpha,\b}(f)\|_\a^2.
\end{align}

\medskip
%By the same idea for \eqref{HRm}, we can obtain higher order Hardy-Rellich equalities, hence inequalities with logarithmic weights by iteration. For example,
\begin{prop}
\label{prop5.4}
Let $\a \in \R$, $u \in C_c^2(B_1\backslash\{0\})$ or $u \in C_c^2(\R^n \backslash\overline{B_1})$. Then
\begin{align}
\label{HRln2a}
\begin{split}
\left\|\Delta u\right\|_{\a}^2 & = A_\a^2\|u\|_{\a+4}^2 + \frac{D_\a}{4}\left\||\ln r|^{-1}u\right\|_{\a+4}^2 + \frac{9}{16}\left\||\ln r|^{-2}u\right\|_{\a+4}^2\\
&\quad +D_\a\left\| G_{\a+2,0} (u)\right\|_{\a+2}^2 + \left\| G_{\a,0} \circ T_{\a+2}(u)\right\|_\a^2 +\frac{1}{4}\left\||\ln r|^{-1} G_{\a+2, -1} (u)\right\|_{\a+2}^2\\
&\quad + 2A_\a\sum_{j=1}^n\|L_j u\|_{\a+2}^2+\|(\Delta-\Delta_r) u\|_\a^2+2\sum_{j=1}^n\|T_\a(L_j u)\|_\a^2.
\end{split}
\end{align}
\end{prop}

\medskip
\noindent
{\sl Proof}. By \eqref{HR32}, there holds
\begin{align*}
%\label{HR33}
\||\ln r|^\b R_{\g, k}(u)\|_{\g}^2=\frac{(2\b -1)^2}{4}\left\||\ln r|^{\b-1} R_{\g+2, k-1}(u)\right\|_{\g+2}^2 + \||\ln r|^\b G_{\g,\beta}(R_{\g+2, k-1}u)\|_{\g}^2.
\end{align*}
The iteration gives then, for all $k \geq 1$, $\g \in \R$,
\begin{align}\label{HR34}
\begin{split}
\|R_{\g, k}(u)\|_{\g}^2 & =\frac{1}{4}\left\||\ln r|^{-1}R_{\g+2, k-1} u\right\|_{\g+2}^2 + \| G_{\g,0}\circ R_{\g+2, k-1}(u)\|_{\g}^2\\
& = \g_{k+1}^2\left\||\ln r|^{-k-1}u\right\|^2_{\g+2k+2} + \g_k^2\left\||\ln r|^{-k}G_{\g+2k,-k}(u)\right\|_{\g+2k}^2\\
& \quad + \sum_{j=0}^{k-1} \g_j^2 \left\||\ln r|^{-j}G_{\g+2j,-j}\circ R_{\g+2j+2, k-j-1}(u)\right\|_{\g+2j}^2.
\end{split}
\end{align}
Here $\g_0 = 1$, and $\g_j = \frac{(2j-1)!!}{2^j}$ for $j \geq 1$.
%As application, recall \eqref{HR13},
%\begin{align*} \|\Delta_r f\|_\a^2=\frac{(n+\a)^2}{4}\|\p_r f\|_{\a+2}^2 + \|T_\a(\p_r f)\|_\a^2.\end{align*}
Combining \eqref{k_a} and \eqref{HR32}, \eqref{HR34} with $k=1$, we get
\begin{align*}
\|T_\a(\p_r u)\|_\a^2 & =\frac{(n-4-\a)^2}{4}\|T_{\a+2}(u)\|_{\a+2}^2 + \|R_{\a, 1}(u)\|_\a^2\\
& = \frac{(n-4-\a)^2}{4} \left[\frac{1}{4}\left\||\ln r|^{-1} u\right\|_{\a+4}^2+\|G_{\a+2,0} (u)\|_{\a+2}^2 \right]\\
&\quad + \frac{9}{16}\left\||\ln r|^{-2} u\right\|_{\a+4}^2+ \frac{1}{4}\left\||\ln r|^{-1} G_{\alpha+2,-1}(u)\right\|_{\a+2}^2+\|G_{\a,0}\circ T_{\a+2} (u)\|_{\a}^2.
\end{align*}
On the other hand, thanks to \eqref{HR32} with $\beta=0$ and \eqref{rHR1a}, there holds
\begin{align*}
%\label{ln3}
\|\p_r u\|_\a^2 =  \frac{1}{4}\||\ln r|^{-1} u\|_{\a+2}^2 + \frac{(n-2-\a)^2}{4}\|u\|_{\a+2}^2 + \|G_{\a, 0}(u)\|_\a^2.
\end{align*}
Then we deduce \eqref{HRln2a} by the above two identities and \eqref{HRm} with $m=1$. \qed

%Clearly, we get Corollary 5.3 and 5.4 in \cite{CM} with $\a=0$.
\medskip
\noindent
{\sl Proof of Theorem \ref{thmHRln}}.
%The inequality \eqref{HR2ln} is a direct consequence of the identity \eqref{HRln2a} and Lemma \ref{lemSphere}, since we need only to have $2A_\a + (n-1) \geq 0$ which is equivalent to $|\a+2| \leq \sqrt{n^2 - 2n + 2}$.
To get \eqref{HR21ln}, we make use of an identity coming from \eqref{HR32} and \eqref{HR34}, that is
\begin{align*}
&\quad \|\D u\|_\a^2 - \frac{(n+\a)^2}{4} \left\|\nabla u\right\|^2_{\a+2}\\& = \frac{(n-4-\a)^2}{4}\|T_{\a+2} (u)\|_{\a+2}^2 +\|R_{\a,1}(u)\|_{\a}^2 +  \left\|(\Delta - \D_r) u\right\|_\a^2\\&\quad +\frac{(n+\a)(n-8-3\a)}{4} \sum_{j=1}^n \left\|L_j u\right\|^2_{\a+2} +2\sum_{j=1}^n\|T_\a(L_j u)\|_\a^2\\
& = \frac{(n-4-\a)^2}{16}\||\ln r|^{-1}u\|_{\a+4}^2 + \frac{9}{16} \||\ln r|^{-2} u\|_{\a+4}^2  + \left\|(\Delta - \D_r) u\right\|_\a^2 \\& \quad + \frac{(n+\a)(n-8-3\a)}{4} \sum_{j=1}^n \left\|L_j u\right\|^2_{\a+2} +\frac{1}{4} \left\|\frac{G_{\a+2,-1} (u)}{\ln|x|}\right\|_{\a+2}^2+\|G_{\a,0}\circ T_{\a+2}(u)\|_{\a}^2 \\&\quad +2\sum_{j=1}^n\|T_\a(L_j u)\|_\a^2 + \frac{(n-4-\a)^2}{4}\|G_{\a+2,0}(u)\|_{\a+2}^2.
\end{align*}
Applying \eqref{Sphere} once again, we conclude \eqref{HR21ln} with $\frac{(n+\a)(n-8-3\a)}{4} + n-1 \ge 0$. \qed

\subsection{Rellich equality and inequalities in Hyperbolic space}
\label{sect-hyper}
%\reset
Consider the Poincar\'e's hyperbolic ball $\H^n$ with $n \ge 2$. Let $\rho=\ln\frac{1+r}{1-r}$ denote the distance from origin to $x\in \H^n$. Let $\L_j$, $\D_{\rho, \H}$ be the sphere derivatives and the radial Laplacian in $\H^n$ defined respectively by \eqref{Ljhyper}, \eqref{Drhyper}. Here we will prove the Rellich identity \eqref{HR2=hyper}, then apply it to establish new Rellich inequalities, especially in $\H^3$ and $\H^4$. In this subsection, $\|\cdot\|$ denotes the $L^2$ norm in $\H^n$, and $\langle\cdot, \cdot\rangle_\H$ means the corresponding inner product.

\medskip
\noindent
{\sl Proof of Theorem \ref{thm6}}. At first, we observe that the Beltrami-Laplace operator $\Delta_\H$ in $\H^n$ satisfies
\begin{align*}
\Delta_\H - \Delta_{\rho, \H} = \frac{1}{({\rm sh}\rho)^2}\Delta_{{\mathbb S}^{n-1}} = \frac{r^2}{({\rm sh}\rho)^2}\sum_{j=1}^nL_j^2 = \frac{(1-r^2)^2}{4}\sum_{j=1}^nL_j^2 = \sum_{j=1}^n\L_j^2.
\end{align*}
%\medskip
%{\sl Proof of Theorem \ref{compahyper}}.
\medskip
Take $u \in C_c^2(\H^n\backslash\{0\})$, by definition,
\begin{equation}\label{4.1}
\|\Delta_\H u\|^2= \|\Delta_{\rho,\H} u\|^2+ \|(\Delta_\H -\Delta_{\rho,\H}) u\|^2 + 2\Big\langle\Delta_{\rho,\H} u, \sum_{j=1}^n \L^2_j u\Big\rangle_\H.
\end{equation}
Using Lemma \ref{lemLj}, $\sum_{1\le j \le n} x_j\L_j = 0$, $\L_j(r) = 0$.
Similar to the Euclidean case, direct calculation yields
\begin{align}\label{4.2}
\Big\langle\Delta_{\rho,\H} u, \sum_{j=1}^n \L^2_j u\Big\rangle_\H = -\sum_{j=1}^n \Big\langle \L_j(\Delta_{\rho,\H} u), \L_j u \Big\rangle_\H.
\end{align}
Thanks to Lemma \ref{lemLj} (ii), there holds
\begin{align*}
\L_j(\p_\rho)-\p_\rho(\L_j) & = \frac{1-r^2}{2}L_j\Big(\frac{1-r^2}{2}\p_r\Big) - \frac{1-r^2}{2}\p_r\Big(\frac{1-r^2}{2}L_j\Big)\\
& = \frac{(1-r^2)^2}{4}\Big[L_j(\p_r) - \p_r(L_j)\Big] + \frac{r(1-r^2)}{2}L_j\\
& = \frac{(1-r^2)(1+r^2)}{4r}L_j\\
& = (\coth\rho)\L_j.
\end{align*}
Hence
\begin{align*}
\L_j(\Delta_{\rho,\H})& =\L_j\Big[({\rm sh}\rho)^{1-n}\p_\rho(({\rm sh}\rho)^{n-1} \p_\rho)\Big]\\
& ={\rm (sh\rho)}^{1-n} \big(\p_\rho + \coth\rho\big)\Big[({\rm sh}\rho)^{n-1}\big(\p_\rho+\coth\rho\big)\Big]\L_j\\
& = \p^2_\rho(\L_j)+  (n+1) \coth\rho \p_\rho(\L_j) + \Big[n(\coth\rho)^2 - ({\rm sh}\rho)^{-2}\Big] \L_j\\
& = \p^2_\rho(\L_j)+  (n+1) \coth\rho \p_\rho(\L_j) + \frac{n-1}{({\rm sh}\rho)^2} \L_j + n \L_j.
\end{align*}
Moreover, we have
%(as $dv_\H = ({\rm sh}\rho)^{n-1} \rho^{1-n}dx$)
\begin{align}
\label{ippH}
\langle \p_\rho f, h\rangle_\H = -\langle \p_\rho h, f\rangle_\H - (n-1)\langle f, (\coth\rho)h\rangle_\H, \quad \forall\; f, h \in C_c^1(\H^n\backslash\{0\}).
\end{align}
Then
\begin{align*}
& \quad \langle \L_j(\Delta_{\rho,\H} u), \L_j u \rangle_\H\\& = \Big\langle \p^2_\rho(\L_j u) + (n+1) \coth\rho \p_\rho(\L_j u) + \frac{n-1}{({\rm sh}\rho)^2} \L_j u + n \L_j u, \L_j u\Big\rangle_\H\\
& = -\|\p_\rho(\L_j u)\|^2 + 2 \langle\p_\rho(\L_j u), \coth\rho\L_j u\rangle_\H + n\|\L_j u\|^2 + (n-1)\Big\|\frac{\L_j u}{{\rm sh}\rho}\Big\|^2\\
& = -\|\p_\rho(\L_j u)\|^2 + \Big\|\frac{\L_j u}{{\rm sh}\rho}\Big\|^2 - (n-1)\|\coth\rho\L_j u\|^2 + n\|\L_j u\|^2 + (n-1)\Big\|\frac{\L_j u}{{\rm sh}\rho}\Big\|^2\\
& = -\|\p_\rho(\L_j u)\|^2 + \Big\|\frac{\L_j u}{{\rm sh}\rho}\Big\|^2 + \|\L_j u\|^2.
\end{align*}
Combining with \eqref{4.1} and \eqref{4.2}, we get the identity \eqref{HR2=hyper}.

\medskip
Furthermore, applying Lemma \ref{lemSphere}, for any $u\in C_c^2(\H^n \backslash \{0\})$, there holds
\begin{align}
\label{D-Drhyper}
\begin{split}
\|(\Delta_\H - \Delta_{\rho, \H})u\|^2 & = \int_0^\infty\int_{\Sp_r^{n-1}} \frac{(1-r^2)^4}{16}\Big|\sum_j L_j^2 u\Big|^2 ({\rm sh}\rho)^{n-1}\rho'(r)d\sigma dr\\
%& = \int_0^\infty \frac{(1-r^2)^4}{16}({\rm sh}\rho)^{n-1}\rho'(r)\int_{\Sp_r^{n-1}} \Big|\sum_j L_j^2 u\Big|^2d\sigma dr\\
& \geq \int_0^\infty \frac{(1-r^2)^4}{16}({\rm sh}\rho)^{n-1}\rho'(r)\frac{n- 1}{r^2}\int_{\Sp_r^{n-1}}\sum_{j= 1}^n |L_j u|^2d\sigma dr\\
& = (n-1)\sum_{j= 1}^n\Big\|\frac{1-r^2}{2r}\L_j u\Big\|^2\\
& = (n-1)\sum_{j= 1}^n\Big\|\frac{\L_j u}{{\rm sh}\rho}\Big\|^2.
\end{split}
\end{align}
 %The result seems to be new in $\R^3$.\footnote{\; The formula holds also in $\H^2$ when $u$ is supported in a small punctured ball around the origin}

\medskip
In addition, similar to \eqref{rHR1a}, the following equality in $\H^n$ holds true, see \cite[Theorem 1.4]{FLLM}. For any $u \in C_c^1(\H^n\backslash\{0\})$,
%, which improves Theorem 2.1 in \cite{BGG}.
%\footnote{This equality could exist in \cite{BGR}}.
\begin{align}
\label{hyp}
\begin{split}
\|\partial_\rho u\|^2 & = \frac{1}{4}\Big\|\frac{u}{\rho}\Big\|^2 + \frac{(n-1)^2}{4}\|u\|^2 + \frac{(n-1)(n-3)}{4}\Big\|\frac{u}{{\rm sh}\rho}\Big\|^2\\
& \quad +\int_{\H^n} \Big|\partial_\rho \Big(\frac{u}{\rho^{\frac{1}{2}}({\rm sh} \rho)^{\frac{1-n}{2}}}\Big)\Big|^2 \rho ({\rm sh} \rho)^{1-n}d v_\H.
\end{split}
\end{align}
Combining \eqref{HR2=hyper}, \eqref{D-Drhyper}, and \eqref{hyp} with $\L_j u$, we get finally
\begin{align*}
%\label{HRhyper1}
\begin{split}
\|\Delta_\H u\|^2& \geq \|\Delta_{\rho,\H} u\|^2 +\frac{1}{2}\sum_{j=1}^n\Big\| \frac{\L_j u}{\rho}\Big\|^2 +\frac{(n+1)(n-3)}{2}\sum_{j=1}^n\|\L_j u\|^2\\
& \quad +\frac{(n+1)(n-3)}{2}\sum_{j=1}^n\Big\|\frac{\L_j u}{{\rm sh}\rho}\Big\|^2,\\
&=\|\Delta_{\rho,\H} u\|^2 +\frac{1}{2}\sum_{j=1}^n\Big\| \frac{\L_j u}{\rho}\Big\|^2 +\frac{(n+1)(n-3)}{2}\sum_{j=1}^n\|{\rm coth}\rho\L_j u\|^2.
\end{split}
\end{align*}
The proof is completed.
\qed

\medskip
\noindent
{\sl Proof of Theorem \ref{prop9}}.
%By \eqref{Rhyper} and Fubini theorem, we only need to prove the same estimate for radial functions.
Firstly, let $u$ be a {\it radial} function in $C_c^2(\H^n\backslash\{0\})$. Denote $v=({\rm sh} \rho)^{\frac{n-1}{2}}u$, by definition of $\Delta_{\rho,\H}$, there holds
\begin{align*}
\Delta_{\rho,\H} u & = ({\rm sh}\rho)^{-\frac{n-1}{2}}\left[\frac{\p^2 v}{\p\rho^2} - \frac{(n-1)(n-3)}{4}(\coth\rho)^2 v - \frac{n-1}{2}v\right] = ({\rm sh}\rho)^{-\frac{n-1}{2}}\big(v'' - Av\big)
\end{align*}
where $ v'' = \p^2_{\rho}v$ and
\begin{align*}
A(\rho) = \frac{(n-1)^2}{4} + \frac{(n-1)(n-3)}{4({\rm sh}\rho)^2}.
\end{align*}
Then $(\Delta_{\rho,\H} u)^2 = ({\rm sh}\rho)^{1-n}\big(v''^2 + A^2v^2 - 2Avv''\big)$, so that
\begin{align*}
\|\Delta_\H u\|^2 = \o_n\int_0^\infty \big(v''^2 + A^2v^2 - 2Avv''\big)d\rho.
\end{align*}
Clearly
\begin{align*}
\int_0^\infty - 2Avv''d\rho = \int_0^\infty \big(2Av'^2 - A''v^2\big)d\rho,
\end{align*}
and
$$A''(\rho) = \frac{(n-1)(n-3)}{4}\Big[\frac{6}{({\rm sh}\rho)^4} + \frac{4}{({\rm sh}\rho)^2}\Big].$$
For the integral of $Av'^2$, define $w = \frac{v}{{\rm sh}\rho}$, then $\frac{v'}{{\rm sh}\rho} = w' + w\coth\rho$. Hence we get
\begin{align*}
\frac{v'^2}{({\rm sh}\rho)^2} & = w'^2 + 2ww'\coth\rho + w^2\coth^2\rho = w'^2 + \big(w^2\coth\rho\big)' + \frac{2v^2}{({\rm sh}\rho)^4} + \frac{v^2}{({\rm sh}\rho)^2}.
\end{align*}
Consequently
\begin{align*}
\int_0^\infty 2Av'^2 d\rho = \frac{(n-1)^2}{2}\int_0^\infty v'^2 d\rho + \frac{(n-1)(n-3)}{2}\int_0^\infty \Big[w'^2 + \frac{2v^2}{({\rm sh}\rho)^4} + \frac{v^2}{({\rm sh}\rho)^2}\Big]d\rho.
\end{align*}
Hence
\begin{align}
\label{newhyp1}
\begin{split}
\int_0^\infty \big(v''^2 + A^2v^2 - 2Avv''\big)d\rho & = \int_0^\infty v''^2 d\rho +  \frac{(n-1)^4}{16}\int_0^\infty v^2 d\rho + \frac{(n-1)^2}{2}\int_0^\infty v'^2 d\rho\\
& \quad +  \Big[\frac{(n-1)^2(n-3)^2}{16} - \frac{(n-1)(n-3)}{2}\Big]\int_0^\infty \frac{v^2}{({\rm sh}\rho)^4}d\rho\\
& \quad +  \Big[\frac{(n-1)^3(n-3)}{8} - \frac{(n-1)(n-3)}{2}\Big]\int_0^\infty \frac{v^2}{({\rm sh}\rho)^2}d\rho\\
& \quad + \frac{(n-1)(n-3)}{2}\int_0^\infty w'^2d\rho.
\end{split}
\end{align}
Applying the classical Hardy inequality over $(0, \infty)$, we see that
\begin{align*}
\int_0^\infty v''^2 d\rho \geq \frac{9}{16}\int_0^\infty \frac{v^2}{\rho^4}d\rho, \quad \int_0^\infty v'^2 d\rho \geq \frac{1}{4}\int_0^\infty \frac{v^2}{\rho^2}d\rho, \quad \int_0^\infty w'^2d\rho \geq \frac{1}{4}\int_0^\infty \frac{v^2}{\rho^2({\rm sh}\rho)^2}d\rho.
\end{align*}
Combining with the identity \eqref{newhyp1}, the verification of \eqref{HRhyper2} is completed for radial function.

\smallskip
%Moreover, as $\rho \leq {\rm sh}\rho$ in $\R_+$, \eqref{HRhyper1} yields that $\|\Delta_\H u\| \ge \|\Delta_{\rho,\H} u\|$ for any $u \in C_c^2(\H^n\backslash\{0\})$ and $n \ge 4$.
By the Fubini theorem, the above estimates work also to get \eqref{HRhyper2} with $\Delta_{\rho, \H} u$, hence with $\Delta_\H u$ and $n \geq 3$, thanks to Theorem \ref{thm6}.

\medskip
 The estimate in $\H^3$ is an obvious consequence of \eqref{HRhyper2}. To get the inequality in $\H^4$, we use
$$\frac{9}{16\rho^4} - \frac{15}{16({\rm sh}\rho)^4} + \frac{3}{8\rho^2({\rm sh}\rho)^2} = \frac{3}{16\rho^4}\Big[3 + \frac{5\rho^2}{({\rm sh}\rho)^2}\Big]\Big[1 - \frac{\rho^2}{({\rm sh}\rho)^2}\Big] >0, \quad\forall\rho > 0,$$
since $\rho^2 \le ({\rm sh}\rho)^2$ in $(0, \infty)$. So we are done.\qed

\bigskip
\noindent
{\bf Acknowledgements.} The authors are partially supported by NSFC (No.~12271164) and Science and Technology Commission of Shanghai Municipality (No.~22DZ2229014).

\end{document}